\documentclass[review]{elsarticle}

\usepackage{lineno,hyperref}
\modulolinenumbers[5]

\journal{Journal of \LaTeX\ Templates}

\usepackage{amssymb}
\usepackage{latexsym}
\usepackage{mathtools}
\usepackage{url}
\usepackage{xcolor}
\definecolor{newcolor}{rgb}{.8,.349,.1}

\usepackage{amsmath}
\usepackage{mathtools}
\usepackage{amsfonts}
\usepackage{amsthm}
\usepackage{bm}
\usepackage{physics}

\newtheorem{lemma} {Lemma}

\theoremstyle{remark}

\theoremstyle{definition}

\usepackage{floatrow}

\graphicspath{ {images/} }

\usepackage[capitalise]{cleveref}
\crefformat{equation}{(#2#1#3)}

\newcommand{\R}{\mathbb{R}}

\newcommand{\mat}[1]{\mathbb{#1}}
\newcommand{\dual}[1]{\tilde{#1}}
\let\abs\relax
\newcommand{\abs}[1]{|{#1}|}
\newcommand{\bdual}[1]{\tilde{\bm{#1}}}
\newcommand{\mc}{\mathcal}

\newcommand{\restr}[2]{#1_{\mkern 1mu \vrule height 1ex\mkern 2mu #2}}

\newcommand{\wrestrs}[3]{\restr{{{}#1_#2}}{#3}}
\newcommand{\wrestr}[2]{\restr{{{}#1}}{#2}}

\newcommand{\f}{\mc{F}}
\newcommand{\e}{\mc{E}}

\newcommand{\df}{{\tilde{\mc{F}}}}
\newcommand{\de}{{\tilde{\mc{E}}}}

\newcommand{\DC}{{\dual{c}}}
\newcommand{\C}{c}

\newcommand{\E}{\restr{E}{\C}}
\newcommand{\F}{\restr{F}{\C}}

\newcommand{\ED}{\restr{E}{\DC}}
\newcommand{\FD}{\restr{F}{\DC}}

\newcommand{\DED}{\restr{\dual{E}}{\DC}}
\newcommand{\DFD}{\restr{\dual{F}}{\DC}}

\newcommand{\mE}{\restr{\mat{E}}{\C}}
\newcommand{\mF}{\restr{\mat{F}}{\C}}

\newcommand{\mDE}{\restr{\dual{\mat E}}{\C}}
\newcommand{\mDF}{\restr{\dual{\mat F}}{\C}}

\newcommand{\mED}{\restr{\mat{E}}{\DC}}
\newcommand{\mFD}{\restr{\mat{F}}{\DC}}

\newcommand{\mDED}{\restr{\dual{\mat E}}{\DC}}
\newcommand{\mDFD}{\restr{\dual{\mat F}}{\DC}}

\newcommand{\dofE}[1]{\restr{{#1}^{\mc{E}}}{\C}}
\newcommand{\dofF}[1]{\restr{{#1}^{\mc{F}}}{\C}}
\newcommand{\dofDE}[1]{\restr{{#1}^{\tilde{\mc{E}}}}{\C}}
\newcommand{\dofDF}[1]{\restr{{#1}^{\tilde{\mc{F}}}}{\C}}

\newcommand{\dofED}[1]{\restr{{#1}^{\mc{E}}}{\DC}}
\newcommand{\dofFD}[1]{\restr{{#1}^{\mc{F}}}{\DC}}
\newcommand{\dofDED}[1]{\restr{{#1}^{\tilde{\mc{E}}}}{\DC}}
\newcommand{\dofDFD}[1]{\restr{{#1}^{\tilde{\mc{F}}}}{\DC}}

\newcommand{\GdofE}[1]{{#1}^{\, \e}}
\newcommand{\GdofF}[1]{{#1}^{\, \f}}
\newcommand{\GdofDE}[1]{{#1}^{\, \de}}
\newcommand{\GdofDF}[1]{{#1}^{\, \df}}

\newcommand{\GdofN}[1]{{#1}^{\, \mc N}}
\newcommand{\GdofDN}[1]{{#1}^{\dual{\mc N}}}

\newcommand{\MF}{\restr{\mat M^{\mc F}}{\C}}
\newcommand{\ME}{\restr{\mat M^{\mc E}}{\C}}
\newcommand{\MDF}{\restr{\mat M^{\tilde{\mc F}}}{\DC}}
\newcommand{\MDE}{\restr{\mat M^{\tilde{\mc E}}}{\DC}}

\newcommand{\MFD}{\restr{\mat M^{\mc F}}{\DC}}
\newcommand{\MED}{\restr{\mat M^{\mc E}}{\DC}}

\newcommand{\GMF}{\mat M^{\mc F}}
\newcommand{\GME}{\mat M^{\mc E}}
\newcommand{\GMDF}{\mat M^{\tilde{\mc F}}}
\newcommand{\GMDE}{\mat M^{\tilde{\mc E}}}
\newcommand{\material}{\mat {K}}

\newcommand{\dualG}{{\mat D^T}}
\newcommand{\dualC}{\mat C^T}

\newcommand{\mDFDR}{\restr{\dual{\mat F}}{\restr{\DC}{\C}}}
\newcommand{\mDEDR}{\restr{\dual{\mat E}}{\restr{\DC}{\C}}}

\newcommand{\mSR}{\restr{\mat S}{\restr{\DC}{\C}}}
\newcommand{\mLR}{\restr{\mat L}{\restr{\DC}{\C}}}

\newcommand{\MFR}{\restr{\mat M^{\mc F}}{\restr{\DC}{\C}}}
\newcommand{\MER}{\restr{\mat M^{\mc E}}{\restr{\DC}{\C}}}

\usepackage [autostyle, english = american]{csquotes}
\MakeOuterQuote{"}

\begin{document}


\begin{frontmatter}

\title{
Explicit geometric construction of \emph{sparse inverse mass matrices} for arbitrary tetrahedral grids}

\author[1]{Silvano Pitassi\corref{cor1}}
\ead{pitassi.silvano@spes.uniud.it}
\cortext[cor1]{Corresponding author:
  Tel.: +039-0432-558037;}
\author[1]{Francesco Trevisan}
\ead{francesco.trevisan@uniud.it}
\author[1]{Ruben Specogna}
\ead{ruben.specogna@uniud.it}

\address[1]{University of Udine, Polytechnic Department of Engineering and Architecture, EMCLab, via delle scienze 206, 33100 Udine, Italy}


\begin{abstract}
The geometric reinterpretation of the Finite Element Method (FEM) shows that Raviart--Thomas and N\'{e}d\'{e}lec mass matrices map from degrees of freedoms (DoFs) attached to geometric elements of a tetrahedral grid to DoFs attached to the barycentric dual grid. The algebraic inverses of the mass matrices map DoFs attached to the barycentric dual grid back to DoFs attached to the corresponding primal tetrahedral grid, but they are of limited practical use since they are dense.

In this paper we present a new geometric construction of \emph{sparse inverse mass matrices} for arbitrary tetrahedral grids and possibly inhomogeneous and anisotropic materials, debunking the conventional wisdom that the barycentric dual grid prohibits a sparse representation for inverse mass matrices. In particular, we provide a unified framework for the construction of both edge and face mass matrices and their sparse inverses. Such a unifying principle relies on novel geometric reconstruction formulas, from which, according to a well-established design strategy, local mass matrices are constructed as the sum of a consistent and a stabilization term. 
A major difference with the approaches proposed so far is that the consistent term is defined geometrically and explicitly, that is, without the necessity of computing the inverses of local matrices. This provides a sensible speedup and an easier implementation.
%
%
%
We use these new sparse inverse mass matrices to discretize a three-dimensional Poisson problem, providing the comparison between the results obtained by various formulations on a benchmark problem with analytical solution.
\end{abstract}


\end{frontmatter}

\section{Introduction}

Many physics-compatible discretization methods depend on a primal-dual mesh structure and the related staggered positioning of the physical variables \cite{tontiold}, often referred to as \emph{degrees of freedom} (DoFs).
In particular, most physics-compatible discretizations for the Maxwell's equations introduced so far, like the Yee-scheme \cite{Yee1966}, the Cell Method (CM) \cite{tontiold}, \cite{Tonti2001}, the Finite Integration Technique \cite{Clemens2001}, the Discrete Geometric Approach (DGA) \cite{cmame}, \cite{jcp}, the Compatible Discrete Operators (CDO) \cite{bonellem2an} are based on a pair of interlocked primal and dual grids. In other physics-compatible methods like the Mimetic Finite Difference method \cite{Lipnikov2014} and the Finite Element Method (FEM) based on Whitney edge and face basis functions (i.e. the N\'{e}d\'{e}lec curl-conforming basis functions and the Raviart--Thomas div-conforming basis functions) \cite{Bossavit1988}, \cite{bossavitbook}, a dual grid is not explicitly used but an interpretation using it is readily available \cite{howweak}, \cite{bossavitbook}, at least for the lowest order version.

A fundamental principle of compatible numerical schemes is the distinction between \emph{topological equations} and \emph{constitutive equations} \cite{tontiold}, \cite{tontinew}. Topological equations express conservation laws of physical theories. As the name suggests, topological equations are valid under arbitrary homeomorphic transformations of the domain in which they are defined.
Constitutive equations, also called \emph{material equations}, describe the behaviour of material, substance or medium under interest.
The \emph{discrete Hodge star operator} \cite{kettunen}, \cite{bossavitbook} acts as discrete counterpart of the material parameters used to formulate constitutive equations and its form depends on the analytic form of the material parameters and on the geometric properties of the grid which discretizes the computational domain. In particular, the algebraic realization of the Hodge operator is given by the so-called \emph{mass matrices} \cite{howweak}, \cite{kettunen}, \cite{bossavitbook}.

What is usually achieved is the construction of discrete Hodge operators which map DoFs attached to primal faces to DoFs attached to dual edges or DoFs attached to primal edges to DoFs attached to dual faces \cite{kettunen}, \cite{bossavitbook}, \cite{cmame}, \cite{jcp}.
This is exactly what the FEM mass matrices computed with Raviart--Thomas and N\'{e}d\'{e}lec basis functions perform, \cite{howweak}, \cite{bossavitbook}, respectively.

However, we can nonetheless define \emph{inverse discrete Hodge operators} that map DoFs attached to dual geometric elements to DoFs attached to primal geometric elements. Explicitly computing the algebraic inverse of the mass matrix is not considered a viable solution given that such a matrix would be dense, so of questionable usefulness in practice.
Therefore, we are in particular interested in a \emph{sparse realization} of inverse discrete Hodge operators.
%
Inverse Hodge operators play a vital role in many applications. Most notably, inverse Hodge operators enable consistent and explicit schemes to solve time-domain wave propagation problems \cite{Yee1966}, \cite{Codecasa2008}, \cite{teixeira}, \cite{Codecasa2018}. Other applications enabled by the inverse mass matrices comprise the explicit construction of the \emph{codifferential operator}, the \emph{Laplace--de Rham operator} \cite{bell} and compute the discrete Hodge decomposition of discrete fields \cite{bell}, \cite{caltagirone}.


The construction of inverse Hodge operators for tetrahedral grids is not entirely new.
In CM \cite{tontiold}, \cite{Tonti2001}, Discrete Exterior Calculus (DEC) \cite{Hirani2003} and in the cell-centered Finite Volume literature \cite{Eymard2000}, a \emph{Vorono\"{\i} dual grid} based on circumcenters is used, and the resulting mass matrices are diagonal in such a way that their inverses can be easily computed.
It is important to note that, when the material is anisotripic, the mass matrices are in general not diagonal, so a sparse inverse Hodge operator cannot be easily constructed.
Another solution using a Vorono\"{\i} dual grid is presented in \cite{bajaj}.
However, most commonly used mesh generators like NETGEN and GMSH produce tetrahedral grids that are not Delaunay and in this case the Vorono\"{\i} dual grid cannot be defined \cite{bellina}.

The aim of this paper is to extend the construction of sparse inverse Hodge operators on \emph{barycentric dual grids} which can be defined for arbitrary tetrahedral grids. The barycentric dual grids are explicitly used in DGA \cite{jcp} and implicitly in the Finite Element Method (FEM) \cite{howweak} and in the MFD \cite{Lipnikov2014}.
Yet, devising a recipe to construct sparse inverse Hodge operators on a barycentric dual grid appears to be a formidable task \cite{bossavitmh}, \cite{teixeira}, \cite{bell}, \cite{bajaj} given that the conventional wisdom is that the barycentric dual grid ``prohibits a sparse representation for their inverse operators'' \cite{hirani2} and, consequently, only approximate constructions have been proposed \cite{teixeira}.

By using a barycentric dual grid, in \cite{Codecasa2008}, \cite{iceaa}, \cite{Codecasa2018}, inverse mass matrices that map from dual faces to primal edges are constructed by assembling local contributions inside dual cells and then computing the algebraic inverse of the resulting local matrices. Yet, we note that the time needed to compute all local inverses is not negligible, because the rank of the local matrices is twenty or more.
Simultaneously with this paper, in \cite{beltman} similar ideas are developed by using the differential forms formalism. Although the differential form approach is more general, it is much more complicated to present. In addition, the theory is developed for homogeneous materials only.

In this paper we provide a unified framework for the construction of both edge and face mass matrices and their sparse inverses.
The unifying principle relies on the construction of discrete Hodge operators starting from \emph{reconstruction formulas} (that provide a representation of the discrete variables at the continuous level) defined by geometric elements of primal and dual grids. Starting from these reconstruction formulas, local mass matrix is constructed according to a well-established design strategy, as the sum of a consistent and a stabilization term \cite{Lipnikov2014}, \cite{bonellecad}. While the consistent term is fixed, the stabilization term depends on some user-dependent parameters, which for instance, offer the possibility of controlling the eingenvalues of the matrices. 

The construction of the inverse mass matrix that maps from dual edges of a barycentric dual grid to primal faces of an arbitrary tetrahedral grid and possibly inhomogeneous materials proposed in this paper is entirely new.
A key difference with respect to the approach proposed in \cite{Codecasa2008}, \cite{iceaa}, \cite{Codecasa2018} in the construction of the other inverse mass matrix that maps from dual faces of the barycentric dual grid to primal edges of the simplicial grid is the \emph{explicit} and geometric construction of the consistent term, without the necessity of computing local inverses.
In addition, the recipe to construct the stabilization term adds more flexibility in the construction of the mass matrices.

In general, the major contribution of our construction is not only to have avoided the computation of the local inverses, but rather the fact that it results in symmetric expressions for both types of local mass matrices, thus highlighting the duality relationship between the pair of grids.
A key requirement of our setting would be that the material parameters are constant on dual cells associated with the nodes of the primal simplicial grid.
Instead, to deal with the general case of material parameters that are constant on the grid cells, but arbitrary discontinuous across the interfaces between the cells, two different approaches are proposed. The first one is a weighted averaging technique classically used in Finite Volumes Methods \cite{Eymard2000}. The second one is based on an extension of the approach proposed \cite{Codecasa2008},\cite{Codecasa2018}.

The paper is organized as follows.
In \cref{compatible}, we introduce the geometric elements from which general tetrahedral grids and their barycentric dual are constructed.
In addition, we present the basic building blocks of low-order compatible numerical schemes that will be used in the work.
In \cref{primal_hodge}, we detail the construction of mass matrices.
In \cref{dual_hodge}, following the same design principles of the previous section, we detail the construction of inverses of mass matrices.
In addition, we describe how to deal with problems having inhomogeneous materials.
In \cref{numerical_results}, we present an application of the newly derived inverse mass matrices to a Poisson problem and we provide a comparison among different formulations that can be used to solve the same problem on the same mesh.
Finally, in \Cref{conclusions}, the conclusions are drawn.

\section{Compatible discretization}
\label{compatible}
A detailed presentation of compatible discretizations can be found in \cite{Bochev2007}. In what follows, we only present the main ideas that will be used throughout the paper.

\subsection{Geometry of primal and dual grid}
\label{grid}
\begin{figure}[!h]
	\centering
	\includegraphics[scale=0.4]{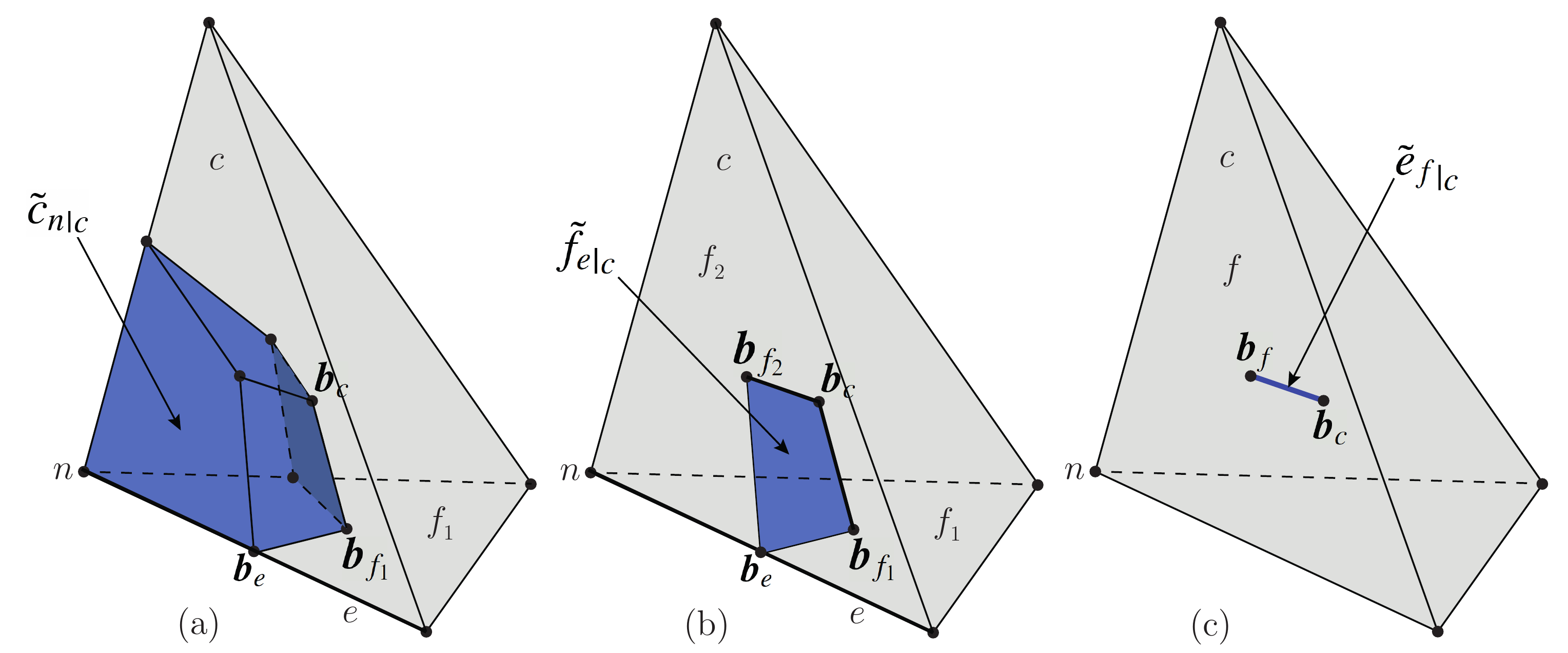}
	\caption{Primal and dual geometric elements of a tetrahedron $\C \in C$. (a) Dual cell. (b) Dual face. (c) Dual edge.}
	\label{cell}
\end{figure}

Let us consider a subdivision of a region $\Omega$ of $\R^3$ into a primal tetrahedral grid $K=(N,E,F,C)$, where the sets $N,E,F$ and $C$ contain the grid \textit{nodes}, \textit{edges}, \textit{faces} and \textit{cells} (or \textit{tetrahedra}), respectively. We denote a single \textit{node} as $n$, an \textit{edge} by $e$, a \textit{face} by $f$ and a \textit{cell} by $\C$, respectively, see \cref{cell}.
The geometric elements of the primal grid are provided with an inner orientation \cite{tontiold}, \cite{tontinew}.

Let $X$ be any set among $N,E,F,$ or $C$. If $r$ is a geometric element of $K$, we denote by $X(r)$ the subset defined by
\begin{align}
X(r) \coloneqq \{x \in X \mid r \subset  x \},
\label{first}
\end{align}
if \cref{first} is not void, or otherwise,
\begin{align}
X(r) \coloneqq \{x \in X \mid x \subset  r \}.
\end{align}
For instance, $C(e) = \{ \C \in C \mid e \subset  \C \}$ is the cluster of cells of $C$ that contain the edge $e$ and $E(\C) = \{e \in E \mid e \subset  c\}$ collects the edges of $c \in C$.

Interlocked with the primal grid $K$, a \textit{barycentric dual grid} \cite{munkres}, \cite{tontiold}, \cite{tontinew} $\dual K = (\dual N,\dual E,\dual F,\dual C)$ is introduced, where the sets $\dual N,\dual E,\dual F$ and $\dual C$ contain \emph{dual nodes}, \emph{dual edges}, \emph{dual faces} and \emph{dual cells}, respectively. Each geometric entity of the dual grid is in one-to-one correspondence (duality pairing) with a geometric element of the primal grid and it is constructed by means of the \emph{barycentric subdivision} \cite{munkres} of the primal grid, see \cref{cell}. With symbol "$\mathtt{\sim}$" we denote geometric elements of the dual grid, thus, we denote a single dual node as $\dual n$, a dual edge by $\dual e$, a dual face by $\dual f$ and a dual cell by $\DC$, respectively. With a subscript we indicate the corresponding (unique) geometric element of the primal grid, thus, the dual of a primal node $n$ is a dual cell $\DC_n$, the dual of a primal edge $e$ is a dual face $\dual f_e$, the dual of a primal face $f$ is a dual edge $\dual e_f$ and the dual of a primal cell $\C$ is the dual node denoted as $\dual n_\C$.

To describe the geometric elements of the pair of grids, we introduce a Cartesian system of coordinates with specified origin and we denote by $\bm x = (x_1,x_2,x_3)^T$ the coordinates of its generic point.
Let us first describe the construction of the \emph{restriction of dual geometric elements to a single cell} $\C$, as the one pictured in \cref{cell}.
The geometric construction of the dual nodes, dual edges, dual faces and dual cells is based on the barycentric subdivision as follows.
To begin with, let us define the barycenters of geometric elements. The barycenters of an edge $e$, face $f$ and of a cell $c$ are defined, respectively, as follows
\begin{align}
&\bm b_e = \frac{1}{\abs{e}} \int_e \bm x\, dl \, ,
&& \bm b_f = \frac{1}{\abs{f}} \int_f \bm x \, dS \, ,
&& \bm b_\C = \frac{1}{\abs{\C}} \int_\C \bm x \, dV \, ,
\end{align}
where $\abs{e},\abs{f}$ and $\abs{c}$ denote the length of edge $e$, the area of face $f$ and the volume of cell $\C$, respectively.

A dual node $\bdual n_\C$ is the barycenter $\bm b_\C$ of $\C$.
A dual edge $\wrestrs{\dual e}{f}{\C}$ is a segment which joins $\bdual n_\C$ with the barycenter $\bm b_f$ of a primal face $f$.
A dual face $\wrestrs{\tilde f}{e}{\C}$ is a quadrilateral plane surface whose vertices are $\bdual n$, the barycenters of the pair of primal faces having a primal edge $e$ in common and the barycenter $\bm b_e$ of $e$.
A dual cell $\wrestrs{\DC}{n}{\C}$ is a region whose boundary is made of a disjoint union of dual faces, corresponding to edges of $\C$ having the node $n$ in common, and portions of faces of $\C$ having the node $n$ in common.
Dual edge $\dual e_f$, dual face $\dual f_e$ and dual cell $\DC_n$ are endowed with outer orientation \cite{tontiold}, \cite{tontinew}, in such a way that each of the pairs $(e, \dual f_e), (f,\dual e_f), (n, \DC_n)$ are oriented according to the right-hand rule.

To each of the following geometric elements $e,f,\wrestr{\dual f}{\C}, \wrestr{\dual e}{\C}$ of the primal or of the dual grid, we associate their corresponding vectors $\bm e, \bm f, \wrestr{\bdual f}{\C}, \wrestr{\bdual e}{\C}$ respectively. Each of these vectors, will be represented with a column array of its cartesian components.
Vector $\bm e$ is the edge vector associated with edge $e$; for example $\bm e = (e_1,e_2,e_3)^T=\bm n_i - \bm n_j, i \neq j$, where $\bm n_i$ and $\bm n_j$ are the coordinates of the boundary nodes of $\bm e$. $\bm t_e$ represents the unit vector parallel to and oriented as $\bm e$ in such a way that $\bm e=\abs{e}\bm t_e$.
Vector $\bm f$ is the face vector associated with face $f$ defined as $\bm f = \abs{f}\bm n_f$, where $\bm n_f$ is the unit vector normal to and oriented as $f$.
In a similar way, vector $\wrestrs{\bdual e}{f}{\C}$ is the edge vector associated with the dual edge $\wrestrs{\dual e}{f}{\C}$; for instance, $\wrestrs{\bdual e}{f}{\C} = \bm b_f - \bdual n_\C$, see \cref{cell}b.
Vector $\wrestrs{\bdual f}{e}{\C}$ is the face vector associated with the dual face $\wrestrs{\dual f}{e}{\C}$.
Face vector $\wrestrs{\bdual f}{e}{\C}$ is equal to $\wrestrs{\bdual f}{e}{\C} =  \frac{1}{2}\big(\wrestrs{\bdual e}{{f_i}}{\C} \times (\bm b_e - \bdual n_\C) - \wrestrs{\bdual e}{{f_j}}{\C} \times (\bm b_e - \bdual n_\C) \big)$ with $f_i,f_j$ two faces such that $e = f_i \cap f_j$, with $i=i(e)$ and $j=j(e)$, and in such a way that the expression induces the correct orientation on $\wrestrs{\bdual f}{e}{\C}$, see \cref{cell}b.

Now, starting from the definition of dual geometric elements restricted to a single cell $\C$, we introduce dual edges, dual faces and dual cells of the barycentric dual grid.
We define a dual edge $\dual e_f$ as the piecewise segment made of the union of the two segments $ \wrestrs{\dual e}{f}{\C_1}$ and $\wrestrs{\dual e}{f}{\C_2}$, with $\C_1,\C_2 \in C(f)$.
The corresponding vector $\bdual e_f$ is defined as follows $\bdual e_f \coloneqq \sum_{\C \in C(f)} \wrestrs{\bdual e}{f}{\C}$.
A dual face $\dual f_e$ is defined as the polyhedral surface made of the union of all $\wrestrs{\dual f}{e}{\C}$, with $\C \in C(e)$.
The corresponding dual face vector $\bdual f_e$ is defined as follows $\bdual f_e \coloneqq \sum_{\C \in C(e)} \wrestrs{\bdual f}{e}{\C}$.
A dual cell $\DC_n$ is defined as the region made of the union of all $\wrestrs{\DC}{n}{\C}$, with $\C \in C(n)$, that is $\DC_n \coloneqq \bigcup_{\C \in C(n)} \wrestrs{\DC}{n}{\C}$.
Note that the boundary of each dual cell decomposes into dual edges and dual faces.

Let $X$ be any set among $E,F$ and let $\dual X$ be any set among $\dual E,\dual F$.
Given a non-empty region $\Omega' \subset \Omega$ of $\R^3$, we define
\begin{align}
&\wrestr{X}{\Omega'} \coloneqq \{ x \cap \Omega', x \in X\},
&\wrestr{\dual X}{\Omega'}\coloneqq \{ \dual x \cap \Omega', \dual x \in X\},
\label{restriction_primal}
\end{align}
as the subsets of geometric elements of $X$ and $\dual X$ restricted to $\Omega'$, respectively.
In the sequel we consider the cases where $\Omega'$ is one among a cell $\C \in C$, a dual cell $\DC \in \dual C$, or the restriction of the dual cell $\DC$ to the cell $\C$, $\wrestr{\DC}{\C}$.
As an instance, $\ED:=\{e \cap \DC, e \in E\}$ contains the restriction of primal edges to the dual cell $\DC$.
Observe that $\E,\F$ coincide with the subsets of all edges and faces of $\C$.
Similarly, $\DED,\DFD$ coincide with the subsets of all dual edges and dual faces of $\DC$. Since dual cells are constructed according to the barycentric subdivision, the length of each edge $e' \in \ED$ is half of the length of a primal edge. Moreover, the area of each face $f' \in \FD$ is a third of the area of a primal face \cite{Berger}.

Finally, we introduce the corresponding matrices $\wrestr{\mat X}{\Omega'}, \wrestr{\dual{\mat X}}{\Omega'}$ whose rows collect geometric vectors associated with elements of the sets in \cref{restriction_primal}.
As an instance, rows of $\mED$ contain edge vectors associated with edges in $\ED$.

\subsection{Degrees of Freedom and reconstruction operators}
\label{dofs}
The physical variables, often referred to as \emph{degrees of freedom} (DoFs), involved in a discrete formulation of a computational problem, have a precise association with the geometric elements of the pair of staggered primal-dual grids \cite{tontiold}, \cite{tontinew}. DoFs are obtained by evaluating the physical vector fields on the geometric elements of the grid, and are collected in finite dimensional arrays of real numbers endowed with a natural vector space structure. The operation of translating a vector field into the corresponding DoFs is performed by means of integration using the so-called \emph{de Rham maps} \cite{kettunen}, \cite{bossavitbook}.
We focus on DoFs attached to edges and faces of the primal and dual grid. We denote the vector space of edge DoFs, face DoFs, dual edge DoFs and dual face DoFs by $\e$, $\f$, $\de$, and $\df$, respectively.
Given a vector field $\bm u$, we denote with $\GdofE{\bm u}$, $\GdofF{\bm u}$, $\GdofDE{\bm u}$, and $\GdofDF{\bm u}$ the DoFs of $\bm u$ in $\e$, $\f$, $\de$, and $\df$, respectively.

We introduce the following specialized definitions of the de Rham maps, which follows as a special case of the general definition applied to the vector subspace of \emph{constant vector fields}.


Let us consider a constant vector field $\bm u$.
Given a non-empty region $\Omega' \subset \Omega$, we define the corresponding subset of DoFs of $\bm u$ \emph{restricted to the region $\Omega'$} as follows
\begin{align}
\restr{{\bm u}^{\e}}{\Omega'} =  \restr{\mat{E}}{\Omega'}  \bm u,\;
\restr{{\bm u}^{\f}}{\Omega'} =  \restr{\mat{F}}{\Omega'}  \bm u,\;
\restr{{\bm u}^{\de}}{\Omega'} = \restr{\dual{\mat E}}{\Omega'} \bm u,\;
\restr{{\bm u}^{\df}}{\Omega'} = \restr{\dual{\mat F}}{\Omega'} \bm u.
\label{local_DOF_primal}
\end{align}

A \emph{reconstruction operator} is defined as a formal inverse of the de Rham map restricted to the vector subspace of constant vector fields. Thus, reconstruction operators map DoFs to a constant vector field defined on a cell. In the sequel we will always refer to local reconstruction operators in such a way that the domain of definition of the reconstructed vector field is a given cell.

\section{Mass matrices}
\label{primal_hodge}
In compatible numerical methods, it is always necessary to map DoFs attached to geometric elements of the primal grid into DoFs attached to geometric elements of the dual grid, or viceversa \cite{tontiold}, \cite{howweak}, \cite{tontinew}.
This process of converting one type of DoFs to another is performed by the so-called \emph{discrete Hodge star operator} \cite{kettunen}.
This operator is required to satisfy a set of properties.
A first requirement is \emph{consistency}, which will be specified in the next subsections, and \emph{symmetry}.
In addition, we require \emph{positive definiteness}, which assures the stability of the resulting numerical scheme.

One possible method for explicitly converting one type of DoFs into another is to reconstruct the polynomial vector field in each cell starting from DoFs attached to a set of geometric elements, and then project this vector field onto the corresponding geometric elements related by duality.
We will follow this principle to design both types of local mass matrices, which give a discrete realization of the corresponding discrete Hodge operators.


In this section we show how to construct consistent mass matrices which map DoFs attached to primal geometric elements to DoFs attached to dual geometric elements.
Novel proofs are given based on Stokes' theorem.
\subsection{Primal cell reconstruction}
Let us consider a pair of vector fields $\bm u, \bm w : \mathbb R^3 \rightarrow \mathbb R^3$ along with a scalar field $w: \mathbb R^3 \rightarrow \mathbb R$ defined on a given cell $\C$. The following well-known integration by parts formulas hold
\begin{equation}
\int_{\C} \bm u \cdot \grad w\, dV=
-\int_{\C} \div \bm u \, w\, dV  +
\int_{\partial \C} \bm u \cdot \bm n \, w \, dS,
\label{stokes1}
\end{equation}

\begin{equation}
\int_{\C} \bm u \cdot \curl{\bm w}\, dV=
\int_{\C} \curl \bm u \cdot \bm w\, dV +
\int_{\partial \C} \bm u \times \bm n \cdot \bm w \, dS.
\label{stokes2}
\end{equation}

\begin{lemma}[Reconstruction formulas]
\label{reconstruction_primal}
The following tensor identities hold
\begin{align}
&\sum_{f \in F(\C)} \wrestrs{\bdual e}{f}{\C} \otimes \bm f = \abs{\C} \, \mat I_3,\label{face_magic}\\
&\sum_{e \in E(\C)} \wrestrs{\bdual f}{e}{\C} \otimes \bm e = \abs{\C} \, \mat I_3.\label{edge_magic}
\end{align}
\begin{proof}
Let $\bm u, \bm w \in \R^3$.
Thanks to \cref{stokes1}, we have
\begin{align}
\begin{split}
\abs{\C} \, \bm u \cdot \bm w \,
&=\int_\C \bm u \cdot \bm w \, dV
=\int_\C \bm u \cdot \grad(\bm w \cdot (\bm x - \bdual n_\C)) \, dV\\
&=\int_{\partial \C} (\bm u \cdot \bm n) \, (\bm w \cdot (\bm x -\bdual n_\C)) \, dS\\
&=\sum_{f \in F(\C)} (\bm u \cdot \bm f) \, ((\bm b_f - \bdual n_\C) \cdot \bm w).
\label{argument}
\end{split}
\end{align}
By using the definition of $\wrestrs{\bdual e}{f}{\C}$, it follows that
\begin{align}
\abs{\C} \, \bm u \cdot \bm w
=\sum_{f \in F(\C)} (\bm u \cdot \bm f) \, (\wrestrs{\bdual e}{f}{\C} \cdot \bm w).
\end{align}
Since $\bm w$ is arbitrary, \cref{face_magic} is proved.

Thanks to \cref{stokes2}, we have
\begin{align}
\begin{split}
2\abs{\C} \, \bm u \cdot \bm w
&=2\int_\C \bm u \cdot \bm w \, dV
=\int_\C \bm u \cdot \curl \big(\bm w \times (\bm x - \bdual n_\C)\big) \, dV\\
&=\int_{\partial \C} (\bm u \times \bm n) \cdot \big(\bm w \times (\bm x - \bdual n_\C)\big) \, dS\\
&=\sum_{f \in F(\C)} \int_{f}  (\bm u \times \bm n_f) \cdot \big(\bm w \times (\bm x - \bdual n_\C)\big) \, dS \\
&=\sum_{f \in F(\C)} \abs{f} \,(\bm u \times \bm n_f) \cdot (\bm w \times \wrestrs{\bdual e}{f}{\C}).
\label{temp}
\end{split}
\end{align}

Now, by applying the same argument used in \cref{argument} to the vector $\bm u \times \bm n_f$ restricted to a face $f$, we obtain
\begin{equation}
\abs{f}\,\bm u \times \bm n_f = \sum_{e \in E(f)} (\bm u \cdot \bm e) \, (\bm b_e -\bm p),
\label{2D_formula}
\end{equation}
where $\bm p\in \R^3$ is an arbitrary node.

Apply \cref{2D_formula} to every face $f$ in the last term of \cref{temp}, choosing the same node $\bdual n_\C$ as the arbitrary point involved in the formula.
We obtain
\begin{align}
\begin{split}
\label{temp2}
2\abs{\C} \, \bm u \cdot \bm w
&=\sum_{f \in F(\C)} \big(\sum_{e \in E(f)} (\bm u \cdot \bm e) \, (\bm b_e - \bdual n_\C) \big) \cdot (\bm w \times \wrestrs{\bdual e}{f}{\C})\\
&=\sum_{f \in F(\C)} \big( \wrestrs{\bdual e}{f}{\C} \times \big(\sum_{e \in E(f)} (\bm u \cdot \bm e) \, (\bm b_e - \bdual n_\C) \big) \big) \cdot \bm w\\
&=\sum_{e \in E(\C)} (\bm u \cdot \bm e) \,  ((\wrestrs{\bdual e}{{f_i}}{\C} \times (\bm b_e - \bdual n_\C)- \wrestrs{\bdual e}{{f_j}}{\C} \times (\bm b_e - \bdual n_\C)) \cdot \bm w),
\end{split}
\end{align}
where $f_i, f_j$ are the unique faces of $\C$ such that $e = f_i \cap f_j$, for suitable indices $i=i(e)$ and $j=j(e)$, and oriented so that they induce opposite orientations on edge $e$.
Now, dividing by two both members of the last term in \cref{temp2} and using the definition of $\wrestrs{\bdual f}{e}{\C}$, it follows that
\begin{equation}
\abs{\C} \, \bm u \cdot \bm w = \sum_{e \in E(\C)} (\bm u \cdot \bm e) \, (\wrestrs{\bdual f}{e}{\C} \cdot \bm w).
\end{equation}
Since $\bm w$ is arbitrary, \cref{edge_magic} is proved.
\end{proof}
\end{lemma}
\subsection{Local mass matrices}
\label{local_mass_matrices}

Let us focus on a cell $\C$ where two pairs of constant vector fields are defined, namely $(\bm u, \bm v)$ and $(\bm w, \bm z)$. The two pairs are related by two constitutive relations
\begin{align}
&\bm v = \wrestrs{\material}{1}{\C} \bm u, \label{primal1}
\\
&\bm z = \wrestrs{\material}{2}{\C} \bm w,
\label{primal2}
\end{align}
where $\wrestrs{\material}{1}{\C},\wrestrs{\material}{2}{\C}$ are two symmetric positive definite matrices of order $3$, assumed to be uniform in $\C$.

Now, let us introduce the restriction of DoFs of the vector fields to $\C$.
In particular, we attach DoFs to geometric elements of the primal and dual grid as follows $\dofF{\bm u} = \mF \bm u$, $\dofDE{\bm v} = \mDE \bm v$, $\dofE{\bm w} = \mE \bm w$, $\dofDF{\bm z} = \mDF \bm z$.

The local mass matrix $\MF$ maps DoFs of $\bm u$ attached to faces to DoFs of $\bm v$ attached to dual edges of the barycentric dual grid.
We say that $\MF$ is a \emph{consistent mass matrix} if
\begin{align}
\dofDE{\bm v} = \MF \dofF{\bm u}
\end{align}
holds exactly for any pair of constant vector fields $\bm u, \bm v$ satisfying \cref{primal1}.

Similarly, a local mass matrix $\ME$ maps DoFs of $\bm w$ attached to edges to DoFs of $\bm z$ attached to dual faces of the barycentric dual grid.
We say that $\ME$ is a \emph{consistent mass matrix} if
\begin{align}
\dofDF{\bm z} = \ME \dofE{\bm w}
\end{align}
holds exactly for any pair of constant vector fields $\bm w, \bm z$ satisfying \cref{primal2}.

An efficient recipe to construct consistent and symmetric matrices $\MF,\ME$ combines the geometric identities in \cref{reconstruction_primal} with the uniformity of the vector fields.

By applying \cref{reconstruction_primal} and the definitions of DoFs $\dofDE{\bm v},\dofF{\bm u}$ of $\bm v$ and $\bm u$, we have that
\begin{align}
\dofDE{\bm v}
=\mDE \bm v
=\mDE \wrestrs{\material}{1}{\C} \bm u
=\mDE \wrestrs{\material}{1}{\C} \frac{1}{\abs{\C}}(\mDE^T \mF) \bm u
=\frac{(\mDE \wrestrs{\material}{1}{\C} \mDE^T)}{\abs{\C}} \dofF{\bm u},
\end{align}
and hence, it follows that a symmetric and consistent matrix $\MF$ is given by
\begin{align}
\MF = \frac{\mDE \wrestrs{\material}{1}{\C} \mDE^T}{\abs{\C}}.
\label{local_ME}
\end{align}

Similarly, by applying \cref{reconstruction_primal} and the definitions of DoFs $\dofDF{\bm z},\dofE{\bm w}$ of $\bm z$ and $\bm w$, we have that
\begin{align}
\dofDF{\bm z}
=\mDF \bm z
=\mDF \wrestrs{\material}{2}{\C} \bm w
=\mDF \wrestrs{\material}{2}{\C} \frac{1}{\abs{\C}}(\mDF^T \mE) \bm w
=\frac{(\mDF \wrestrs{\material}{1}{\C} \mDF^T)}{\abs{\C}} \dofE{\bm w},
\end{align}
and hence, it follows that a symmetric and consistent matrix $\ME$ is given by
\begin{align}
\ME = \frac{\mDF \wrestrs{\material}{2}{\C} \mDF^T}{\abs{\C}}.
\label{local_MF}
\end{align}

The matrices $\MF$ and $\ME$, defined in \cref{local_ME} and \cref{local_MF}, are symmetric and consistent but are not positive definite.
To achieve this, the idea, developed in \cref{local_mass_primal}, is to add to the consistent positive definite matrices $\MF,\ME$ a \emph{stability matrix}, which is symmetric and positive semidefinite.
The stability matrix coincides with one proposed in the mimetic literature \cite{Lipnikov2014}.
\begin{lemma}
\label{local_mass_primal}
Let $m$ be the cardinality of $F(\C)$ or $E(\C)$.
Let $\wrestr{\material}{\C}$ be a symmetric and positive definite matrix of order $3$. Let $\alpha=(\alpha_1,\ldots,\alpha_{m-3})\in(\R^+)^{m-3}$ be any $(m-3)$-upla of positive real numbers and let $\mat D_\alpha$ be the diagonal matrix whose diagonal entries are $\alpha_1,\ldots,\alpha_{m-3}$. Denote by $w_{\f}$ and $w_{\e}$ any two orthonormal basis of $\rm{im}(\mF)^\perp$ and $\rm{im}(\mE)^\perp$, respectively. Accordingly, we define $\wrestr{\mat W}{\C}^\f$ and $\wrestr{\mat W}{\C}^\e$  to be the matrices whose columns collect the vectors of the bases $w_{\f}$ and $w_{\e}$. Then the following matrices
\begin{align}
&\MF \coloneqq \frac{1}{\abs{\C}}\mDE \wrestr{\material}{\C} \mDE^T +  \wrestr{\mat W}{\C}^\f \mat D_\alpha {\wrestr{\mat W}{\C}^\f}^T,\\
&\ME \coloneqq \frac{1}{\abs{\C}}\mDF \wrestr{\material}{\C} \mDF^T +  \wrestr{\mat W}{\C}^\e \mat D_\alpha {\wrestr{\mat W}{\C}^\e}^T,
\end{align}
are symmetric, consistent and positive definite.
\begin{proof}
It sufficient to show that $\MF$ is positive definite. Let $\bm z\in\R^m$ be such that $\bm z^T\MF\bm z=0$.
In order to prove that $\MF$ is positive-definite, we have to show that $\bm z=\bm 0$. The condition $\bm z^T \MF \bm z=0$ is equivalent to require that $(\mDE^T \bm z)^T\wrestr{\material}{\C}\, \mDE^T \bm z=0$ and $({\wrestr{\mat W}{\C}^\f}^T {\bm z})^T \mat D_\alpha( {\wrestr{\mat W}{\C}^\f}^T{\bm z}) =0$.
Since $\wrestr{\material}{\C}$ is symmetric and positive definite, and each $\alpha_i$ is positive, the latter condition is in turn equivalent to $\mDE^T \bm z=\bm 0$ and $\bm z\in\mathrm{im}(\mF)$.
As a consequence, $\bm z=\mF \bm y$ for some $\bm y\in\R^3$ and $\bm 0=\mDE^T \bm z=\mDE^T \mF \bm y=\abs{\C} \, \bm y$.
Thus $\bm z=\mF \bm y=\bm 0$, as desired.
\end{proof}
\end{lemma}

In the general case of a grid made of more than one cell, the corresponding global mass matrices $\GMF,\GME$ are obtained by assembling, cell by cell, the contributions from the local matrices $\MF$ and $\ME$, respectively.


\section{Sparse inverse mass matrices}
\label{dual_hodge}
To derive inverse mass matrices we will mimic the reasoning of \cref{primal_hodge}.
In particular, we will derive a local mass matrix starting from reconstruction formulas and subsequently construct a global mass matrix by applying a standard assembling process.
The construction of inverse mass matrices is carried out only for tetrahedral grids through the introduction of specialized geometric identities.
This approach differs from the one described in \cref{primal_hodge}.
Indeed, the same arguments used in the proof of \cref{reconstruction_primal} cannot be applied to the barycentric dual grid since dual cells have non-planar geometric elements.

\subsection{Dual cell reconstruction}


In the following, we give novel proofs of reconstruction formula that provide a constant vector field defined on $\DC$ starting from DoFs attached to dual edges and dual faces.

Let $\C$ be a tetrahedron and let $n$ be one of its nodes.
Let $i \in\{1,2,3\}$. For $e_i \in E(n) \cap E(\C)$, consider the unique face $f_i \in F(n) \cap F(\C)$ to which $e_i$ does not belong. In this way, for every $i \in \{1,2,3\}$, there exists, and is unique, $r_i \in \{-1,1\}$ such that $r_i \, \bm f_i \cdot \bm e_i > 0$.
\begin{lemma}[Tensor identity]
\label{fundamental_lemma}

The following tensor identity holds
\begin{align}
\sum_{i=1}^3 r_i (\bm f_i \otimes \bm e_i)
=3\abs{\C} \, \mat I_3.
\label{fundamental_identity}
\end{align}

\begin{proof}
Let $j \in \{1,2,3\}$ and let $\bm m_j$ be the unit vector of $\bm f_j$.
Let us prove that
\begin{align}
\sum_{i=1}^3 r_i (\bm f_i \otimes \bm e_i)\, \bm m_j =
\sum_{i=1}^3 r_i(\bm e_i \cdot \bm m_j) \, \bm f_i = 3 \abs{\C} \, \bm m_j = 3 \abs{\C} \, \mat I_3 \bm m_j,
\end{align}
from which the claimed identity follows, since the set of vectors $\bm m_j$ are linearly independent.
Observe that if $i\neq j$, then $\bm e_i \cdot \bm m_j = 0$. It follows that
\begin{align}
\sum_{i=1}^3 r_i(\bm e_i \cdot \bm  m_j) \, \bm f_i
=r_j(\bm e_j \cdot \bm m_j) \, \bm f_j
=r_j(\bm e_j \cdot \bm m_j)\abs{f_j} \, \bm m_j
=3 \abs{\C} \, \bm m_j,
\end{align}
where we have used the expression of the volume of a tetrahedron \citep{Berger}.
\end{proof}
\end{lemma}

Now, we define the restriction of the following oriented geometric boundary terms $\wrestrs{s}{{n,e}}{\C}, \wrestrs{l}{{n,f}}{\C}$ to a tetrahedron $\C$, which are in bijection with edges $e \in E(n) \cap E(\C)$ and faces $f \in F(n) \cap F(\C)$, respectively.
\begin{figure}[H]
	\centering
	\includegraphics[scale=0.4]{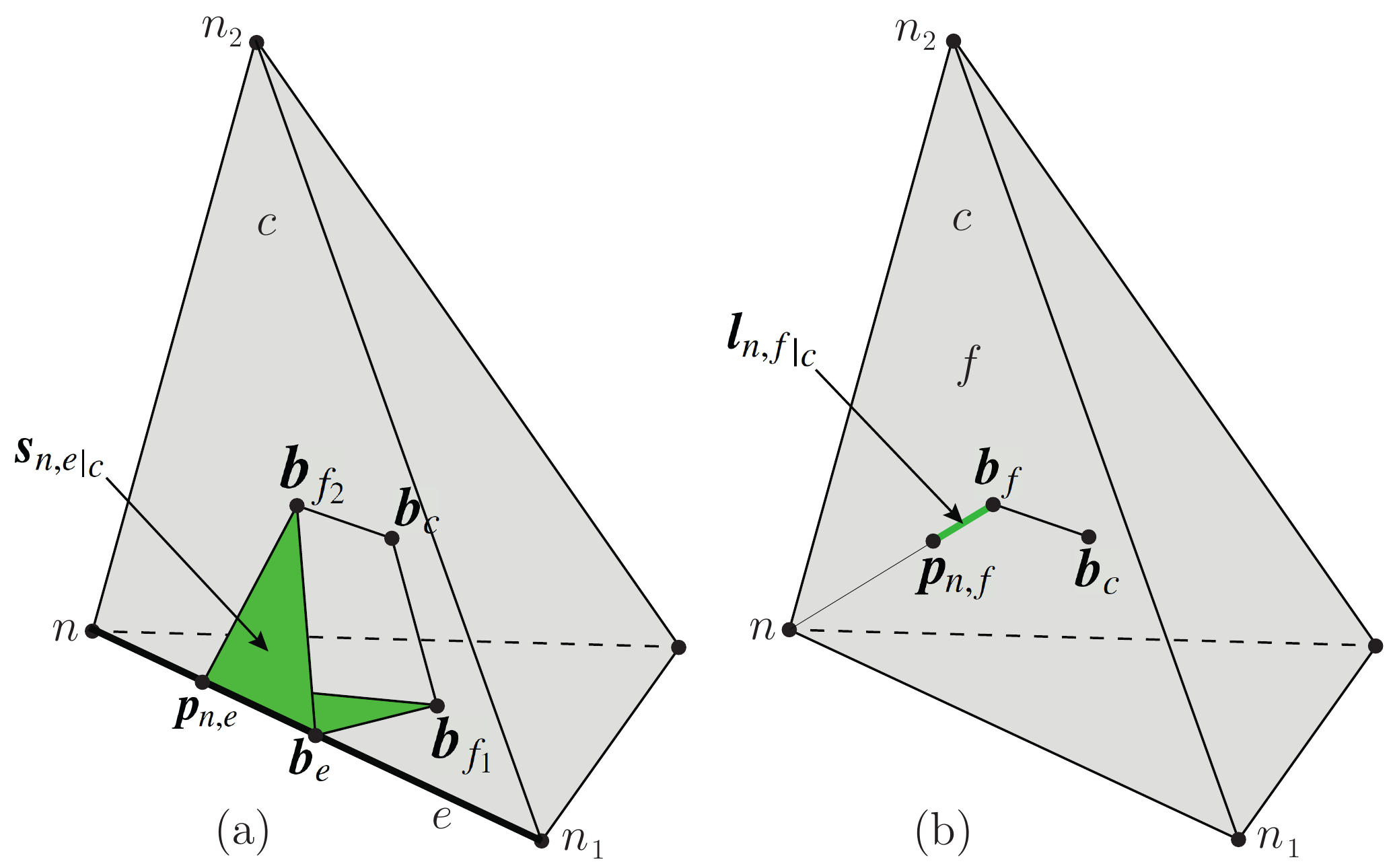}
	\caption{(a) Geometric construction of $\wrestrs{s}{{n,e}}{\C}$. (b) Geometric construction of $\wrestrs{l}{{n,f}}{\C}$.}
	\label{alette}
\end{figure}

We first define $\wrestrs{s}{{n,e}}{\C}$. It decomposes into the union of two triangles, $\wrestrs{t}{{n,e}}{\C}^{(1)}$ and $\wrestrs{t}{{n,e}}{\C}^{(2)}$, associated with the two faces $f_1, f_2 \in F(e)$, see Fig. \ref{alette}(a).
The first has the vertices $\{\bm b_e, \bm b_{f_1}, \bm p_{n,e}\}$, the second the vertices $\{\bm b_e, \bm b_{f_2}, \bm p_{n,e} \}$, with node $\bm p_{n,e}$ defined as follows
\begin{align}
\bm p_{n,e} \coloneqq \frac{3}{4}\bm n +\frac{1}{4}\bm n_1,\, n_1 \in N(e) \setminus \{n\}.
\label{node_e}
\end{align}
By construction, both triangles $\wrestrs{t}{{n,e}}{\C}^{(1)}, \wrestrs{t}{{n,e}}{\C}^{(2)}$ are adjacent to the restriction of the dual face $\wrestrs{\dual f}{e}{\C}$.
Thus, we orient them in such a way that each pair of surfaces $\wrestrs{t}{{n,e}}{\C}^{(1)}, \wrestrs{\dual f}{e}{\C}$ and $\wrestrs{t}{{n,e}}{\C}^{(2)},\wrestrs{\dual f}{e}{\C}$ induce opposite orientations on their common edge \cite{tontinew}.
We define $\wrestrs{s}{{n,e}}{\C}$ to be the surface made of the union of the two triangles $\wrestrs{t}{{n,e}}{\C}^{(1)}, \wrestrs{t}{{n,e}}{\C}^{(2)}$.
We denote with $\wrestrs{\bm s}{{n,e}}{\C}$ the face vector associated with $\wrestrs{s}{{n,e}}{\C}$. It is defined as the sum of the face vectors $\wrestrs{\bm t}{{n,e}}{\C}^{(1)}, \wrestrs{\bm t}{{n,e}}{\C}^{(2)}$, that is $\wrestrs{\bm s}{{n,e}}{\C}  \coloneqq \wrestrs{\bm t}{{n,e}}{\C}^{(1)} +  \wrestrs{\bm t}{{n,e}}{\C}^{(2)}$.
We define $\bm s_{n,e}$ to be $\bm s_{n,e} \coloneqq \sum_{\C \in C(e)} \wrestrs{\bm s}{{n,e}}{\C}$.

Next, we define $\wrestrs{l}{{n,f}}{\C}$ to be the segment joining $\bm b_f$ with node $\bm p_{n,f}$ defined as follows
\begin{align}
\bm p_{n,f} \coloneqq \frac{1}{2} \bm n + \frac{1}{4} \bm n_1 + \frac{1}{4} \bm n_2, \, n_1, n_2 \in N(f) \setminus \{n\},
\label{node_f}
\end{align}
see Fig. \ref{alette}(b).
By construction, $\wrestrs{l}{{n,f}}{\C}$ is adjacent to the restriction of the dual edge $\wrestrs{\dual e}{f}{\C}$.
Thus, we orient $\wrestrs{l}{{n,f}}{\C}$ in such a way the two segments $\wrestrs{l}{{n,f}}{\C}, \wrestrs{\dual e}{f}{\C}$ induce opposite orientations on their common node \cite{tontinew}.
We denote with $\wrestrs{\bm l}{{n,f}}{\C}$ the edge vector associated with $\wrestrs{l}{{n,f}}{\C}$.
We define $\bm l_{n,f}$ to be $\bm l_{n,f} \coloneqq \sum_{\C \in C(f)} \wrestrs{\bm l}{{n,f}}{\C}$.
By mimicking the notation introduced at the end of \cref{grid}, we denote by $\wrestr{\mat S}{\DC}$ and $\wrestr{\mat L}{\DC}$ the matrices whose rows collect vectors $\bm s_{n,e}$ and $\bm l_{n,f}$ with $e \in E(n) $ and $f \in F(n)$, respectively.
Similarly, we denote by $\mSR$ and $\mLR$ the matrices whose rows collect vectors $\wrestrs{\bm s}{{n,e}}{\C}$ and $\wrestrs{\bm l}{{n,f}}{\C}$ with $e \in E(n) \cap E(\C)$ and $f \in F(n) \cap F(\C)$, respectively.


The geometric boundary terms satisfy the following identities
\begin{align}
& \frac{r_i \bm f_i}{6} = \wrestrs{\bdual f}{{e_i}}{\C} + \wrestrs{\bm s}{{n, e_i}}{\C},
\label{first_geom}\\
&\frac{r_i \bm e_i}{4} = \wrestrs{\bdual e}{{f_i}}{\C} + \wrestrs{\bm l}{{n,f_i}}{\C},
\label{second_geom}
\end{align}
which can be proven after trivial algebraic manipulations that we omit.
We point out that if we require $\wrestrs{s}{{n,e}}{\C}, \wrestrs{l}{{n,f}}{\C}$ to satisfy \cref{first_geom}, \cref{second_geom}, the nodes $\bm p_{n,e}, \bm p_{n,f}$ are uniquely identified by \cref{node_e}, \cref{node_f}.

\Cref{first_geom}, \cref{second_geom} are the main geometric identities that together with \cref{fundamental_lemma} allow us to prove the following result.
\begin{lemma}[Reconstruction formulas with boundary terms]
\label{reconsruction_formula_tensor}
We focus attention on a tetrahedron $\C$ and on a dual cell $\DC_n$, dual to a node $n$ of $\C$.
The following tensor identities hold
\begin{align}
\sum_{e \in E(n)\cap E(\C)} \wrestr{\bm e}{\DC} \otimes (\wrestrs{\bdual f}{e}{\C} + \wrestrs{\bm s}{{n,e}}{\C}) = \abs{\wrestr{\DC}{\C}} \, \mat I_3, \label{dual_cell_dual_face}\\
\sum_{f \in F(n)\cap F(\C)} \wrestr{\bm f}{\DC} \otimes (\wrestrs{\bdual e}{f}{\C} + \wrestrs{\bm l}{{n,f}}{\C}) = \abs{\wrestr{\DC}{\C}} \, \mat I_3. \label{dual_cell_dual_edge}
\end{align}
\begin{proof}
Let us apply \cref{fundamental_lemma} to $\C$.
We rewrite \cref{fundamental_identity} as follows
\begin{align}
\frac{1}{12}\sum_{i=1}^3 r_i (\bm e_i \otimes \bm f_i) = \abs{\wrestr{\DC}{\C}} \, \mat I_3,
\end{align}
where we have used the fact that $4\abs{\wrestr{\DC}{\C}} = \abs{\C}$, as a result of barycentric subdivision \cite{Berger}.
Now, using \cref{first_geom} it follows that
\begin{align}
\sum_{i=1}^3 \frac{\bm e_i}{2} \otimes (\wrestrs{\bdual f}{{e_i}}{\C} + \wrestrs{\bm s}{{n,e_i}}{\C})
=\sum_{i=1}^3 \frac{\bm e_i}{2} \otimes \frac{r_i \bm f_i}{6}
=\frac{1}{12}\sum_{i=1}^3 \bm e_i \otimes (r_i \bm f_i)
=\abs{\wrestr{\DC}{\C}} \, \mat I_3.
\end{align}
This proves \cref{dual_cell_dual_face}, since the length of each edge in $\ED$ is half of the length a primal edge.

Similarly, using \cref{second_geom} it follows that
\begin{align}
\sum_{i=1}^3 \frac{\bm f_i}{3} \otimes (\wrestrs{\bdual e}{{f_i}}{\C} + \wrestrs{\bm l}{{n,f_i}}{\C})
=\sum_{i=1}^3 \frac{\bm f_i}{3} \otimes \frac{r_i \bm e_i}{4}
=\frac{1}{12}\sum_{i=1}^3 \bm f_i \otimes (r_i \bm e_i)
=\abs{\wrestr{\DC}{\C}} \, \mat I_3.
\end{align}
This proves \cref{dual_cell_dual_edge}, since the area of each face in $\FD$ is a third of the area of a primal face.
\end{proof}
\end{lemma}

Now we are ready to prove the main result of this section, which is the dual counterpart of \cref{reconstruction_primal}.
This result is remarkable since dual cells have non-planar geometric elements.
\begin{lemma}[Reconstruction formulas]
\label{reconstruction_dual_simplicial}
Let $\DC_n$ be a dual cell.
The following two assertions hold.

(1) If $\DC$ does not intersect with $\partial{\Omega}$, then
\begin{align}
&\sum_{e \in E(n)} \wrestr{\bm e}{\DC} \otimes \bdual f_e  = \abs{\DC} \, \mat I_3,, \label{simplicial_face_internal}\\
&\sum_{f \in F(n)} \wrestr{\bm f}{\DC} \otimes \bdual e_f = \abs{\DC} \, \mat I_3. \label{simplicial_edge_internal}
\end{align}

(2) If $\DC$ does intersect with $\partial{\Omega}$, then
\begin{align}
&\sum_{e \in E(n)} \wrestr{\bm e}{\DC} \otimes (\bdual f_e + \bm s_{n,e}) = \abs{\DC} \, \mat I_3, \label{simplicial_face_external}\\
&\sum_{f \in F(n)} \wrestr{\bm f}{\DC} \otimes (\bdual e_f + \bm l_{n,f}) = \abs{\DC} \, \mat I_3. \label{simplicial_edge_external}
\end{align}

\begin{proof}

Proof of \cref{simplicial_face_internal} and \cref{simplicial_face_external}. Let us apply \cref{dual_cell_dual_face} to every tetrahedron $\C \in C(n)$. By summing over the set $C(n)$ we have
\begin{align}
\sum_{\C \in C(n)} \Bigg( \sum_{e \in E(n) \cap E(\C)} \wrestr{\bm e}{\DC} \otimes (\wrestrs{\bdual f}{e}{\C} + \wrestrs{\bm s}{{n,e}}{\C}) \Bigg)
=\sum_{\C \in C(n)} \abs{\wrestr{\DC}{\C}} \, \mat I_3
=\abs{\DC} \, \mat I_3.
\label{first_sum_face}
\end{align}
We can rewrite \cref{first_sum_face} as a sum over the set $E(n)$ as follows
\begin{align}
\sum_{\C \in C(n)} \Bigg( \sum_{e \in E(n) \cap E(\C)} \wrestr{\bm e}{\DC} \otimes (\wrestrs{\bdual f}{e}{\C} + \wrestrs{\bm s}{{n,e}}{\C}) \Bigg)
=\sum_{e \in E(n)} \wrestr{\bm e}{\DC} \otimes \Bigg(\sum_{\C \in C(e)} \wrestrs{\bdual f}{e}{\C} + \wrestrs{\bm s}{{n,e}}{\C}\Bigg).
\label{second_sum_face}
\end{align}
By using the definition of $\bdual f_e$ and $\bm s_{n,e}$, and combining the two expressions \cref{first_sum_face}, \cref{second_sum_face} we obtain
\begin{align}
\sum_{e \in E(n)} \wrestr{\bm e}{\DC} \otimes (\bdual f_{e} + \bm s_{n,e})
=\abs{\DC} \, \mat I_3.
\end{align}
This proves \cref{simplicial_face_external}.
To prove \cref{simplicial_face_internal}, note that if $\DC$ does not intersect with $\partial{\Omega}$, then
\begin{align}
\bm s_{n,e} = \sum_{e \in C(e)} \wrestrs{\bm s}{{n,e}}{\C} = \bm 0
\end{align}
because each term $\wrestrs{\bm s}{{n,e}}{\C}$ decomposes into the union of two triangles, and each of them appears in the sum exactly two times and with opposite orientation.

Proof of \cref{simplicial_edge_internal} and \cref{simplicial_edge_external}. Let us apply \cref{dual_cell_dual_edge} to every tetrahedron $\C \in C(n)$. By summing over the set $C(n)$ we have
\begin{align}
\sum_{\C \in C(n)} \Bigg( \sum_{f \in F(n) \cap F(\C)} \wrestr{\bm f}{\DC} \otimes (\wrestrs{\bdual e}{f}{\C} + \wrestrs{\bm l}{{n,f}}{\C}) \Bigg)
=\sum_{\C \in C(n)} \abs{\wrestr{\DC}{\C}} \, \mat I_3
=\abs{\DC} \, \mat I_3.
\label{first_sum_face2}
\end{align}
We can rewrite \cref{first_sum_face2} as a sum over the set $F(n)$ as follows
\begin{align}
\sum_{\C \in C(n)} \Bigg( \sum_{f \in F(n) \cap F(\C)} \wrestr{\bm f}{\DC} \otimes (\wrestrs{\bdual e}{f}{\C} + \wrestrs{\bm l}{{n,f}}{\C}) \Bigg)
=\sum_{f \in F(n)} \wrestr{\bm f}{\DC} \otimes \Bigg(\sum_{\C \in C(f)} \wrestrs{\bdual e}{f}{\C} + \wrestrs{\bm l}{{n,f}}{\C}\Bigg).
\label{second_sum_face2}
\end{align}
By using the definition of $\bdual e_f$ and $\bm l_{n,f}$, and combining the two expressions \cref{first_sum_face2}, \cref{second_sum_face2} we obtain
\begin{align}
\sum_{f \in F(n)} \wrestr{\bm f}{\DC} \otimes ( \bdual e_{f} + \bm l_{n,f} )
=\abs{\DC} \, \mat I_3.
\label{purpose_edge}
\end{align}
This proves \cref{simplicial_edge_external}.
To prove \cref{simplicial_edge_internal}, note that if $\DC$ does not intersect with $\partial{\Omega}$, then
\begin{align}
\bm l_{n,f} = \sum_{\C \in C(f)} \wrestrs{\bm l}{{n,f}}{\C} = 0
\end{align}
because each term $\wrestrs{\bm l}{{n,f}}{\C}$ appears in the sum exactly two times and with opposite orientation.
\end{proof}
\end{lemma}
%
%

The results of \cref{reconstruction_dual_simplicial} cannot be extended to arbitrary polyhedral grids.
To see this, it is sufficient to consider a generic 3D polyhedral grid and check that the relations in \cref{reconstruction_dual_simplicial} are not satisfied.
Simple random examples provide readily a counterexample.





\subsection{Local inverse mass matrices}
\label{dual_mass_matrix}

Let us focus on a single dual cell $\DC_n$ where two pairs of constant vector fields are defined, namely $(\bm u, \bm v)$ and $(\bm w, \bm z)$.
For the sake of clarity, we suppose that $\DC$ does not intersect with $\partial \Omega$.
The two pairs are related by two constitutive relations
\begin{align}
\bm v = \wrestrs{\dual \material}{1}{\DC} \bm u,\label{dual1}
\\
\bm z = \wrestrs{\dual \material}{2}{\DC} \bm w,\label{dual2}
\end{align}
where $\wrestrs{\dual \material}{1}{\DC},\wrestrs{\dual \material}{2}{\DC}$ are two symmetric positive definite matrices of order $3$, assumed to be uniform in $\DC$.

Now, let us introduce the restriction of DoFs of the vector fields to $\DC$.
In particular, we attach DoFs to geometric elements of the primal and dual grid as follows $\dofDFD{\bm u} = \mDFD \bm u$, $\dofED{\bm v} = \mED \bm v$, $\dofDED{\bm w} = \mDED \bm w$, $\dofFD{\bm z} = \mFD \bm z$.

The local mass matrix $\MDF$ maps DoFs of $\bm u$ attached to dual faces to DoFs of $\bm v$ attached to edges of the primal grid.
We say that $\MDF$ is a \emph{consistent mass matrix} if
\begin{align}
\dofED{\bm v} = \MDF \dofDFD{\bm u}
\label{cons1}
\end{align}
holds exactly for any pair of constant vector fields $\bm u, \bm v$ satisfying \cref{dual1}.

Similarly, a local mass matrix $\MDE$ maps DoFs of $\bm w$ attached to dual edges to DoFs of $\bm z$ attached to faces of the primal grid.
We say that $\MDE$ is a \emph{consistent mass matrix} if
\begin{align}
\dofFD{\bm z} = \MDE \dofDED{\bm w}
\label{cons2}
\end{align}
holds exactly for any pair of constant vector fields $\bm w, \bm z$ satisfying \cref{dual2}.

An efficient recipe to construct a consistent and symmetric matrices $\MDF,\MDE$ combines the geometric identities in \cref{reconstruction_dual_simplicial} with the uniformity of the vector fields.

By applying \cref{reconstruction_dual_simplicial}, and the definitions of DoFs $\dofED{\bm v},\dofDFD{\bm u}$ of $\bm v$ and $\bm u$, we have that
\begin{align}
\dofED{\bm v}
=\mED \bm v
=\mED \wrestrs{\dual \material}{1}{\DC} \bm u
=\mED \wrestrs{\dual \material}{1}{\DC} \frac{1}{\abs{\DC}}(\mED^T \mDFD) \bm u
=\frac{(\mED \wrestrs{\dual \material}{1}{\DC} \mED^T)}{\abs{\DC}} \dofDFD{\bm u},
\label{step1}
\end{align}
and hence, it follows that a symmetric and consistent matrix $\MDF$ is given by
\begin{align}
\MDF = \frac{\mED \wrestrs{\dual \material}{1}{\DC} \mED^T}{\abs{\DC}}.
\label{local_MDE}
\end{align}

Similarly, by applying \cref{reconstruction_dual_simplicial}, and the definitions of DoFs $\dofFD{\bm z},\dofDED{\bm w}$ of $\bm z$ and $\bm w$, we have that
\begin{align}
\dofFD{\bm z}
=\mFD \bm z
=\mFD \wrestrs{\dual \material}{2}{\DC} \bm w
=\mFD \wrestrs{\dual \material}{2}{\DC} \frac{1}{\abs{\DC}}(\mFD^T \mDED) \bm w
=\frac{(\mFD \wrestrs{\dual \material}{1}{\DC} \mFD^T)}{\abs{\DC}} \dofDED{\bm w},
\label{step2}
\end{align}
and hence, it follows that a symmetric and consistent matrix $\MDE$ is given by
\begin{align}
\MDE = \frac{\mFD \wrestrs{\dual \material}{2}{\DC} \mFD^T}{\abs{\DC}}.
\label{local_MDF}
\end{align}

The matrices $\MDF$ and $\MDE$ defined in \cref{local_MDE}, \cref{local_MDF} are symmetric and consistent but are not positive definite.
To achieve this, we add a stability matrix just as we did for the primal mass matrix.

\begin{lemma}
\label{local_mass_dual}
Let $m$ be the cardinality of $F(n)$ or $E(n)$.
Let $\wrestr{\dual \material}{\DC}$ be a symmetric and positive definite matrix of order $3$. Let $\alpha=(\alpha_1,\ldots,\alpha_{m-3})\in(\R^+)^{m-3}$ be any $(m-3)$-upla of positive real numbers and let $\mat D_\alpha$ be the diagonal matrix whose diagonal entries are $\alpha_1,\ldots,\alpha_{m-3}$. Denote by $w_{\dual{\f}}$ and $w_{\dual{\e}}$ any two orthonormal basis of $\rm{im}(\mDFD + \wrestr{\mat S}{\DC})^\perp$ and $\rm{im}(\mDED + \wrestr{\mat L}{\DC})^\perp$, respectively. Accordingly, we define $\wrestr{\mat W}{\DC}^\df$ and $\wrestr{\mat W}{\DC}^\de$  to be the matrices whose columns collect the vectors of the bases $w_{\dual{\f}}$ and $w_{\dual{\e}}$. Then the following matrices
\begin{align}
&\MDF \coloneqq \frac{1}{\abs{\DC}}\mED \wrestr{\dual \material}{\DC} \mED^T +  \wrestr{\mat W}{\DC}^\df \mat D_\alpha {\wrestr{\mat W}{\DC}^\df}^T, \label{edgedual}\\
&\MDE \coloneqq \frac{1}{\abs{\DC}}\mFD \wrestr{\dual \material}{\DC} \mFD^T +  \wrestr{\mat W}{\DC}^\de \mat D_\alpha {\wrestr{\mat W}{\DC}^\de}^T ,\label{facedual}
\end{align}
are symmetric, consistent and positive definite.
\begin{proof}
\end{proof}
It is a direct consequence of \cref{local_mass_primal}.
\end{lemma}

In the general case of a dual grid made of more than one dual cell, the corresponding global mass matrices $\GMDF,\GMDE$ are obtained by assembling, dual cell by dual cell, the contributions from the local matrices $\MDF$ and $\MDE$, respectively.
However, we need also to take into account the case where $\DC$ does intersect with $\partial{\Omega}$.
In this case, \cref{simplicial_face_external}, \cref{simplicial_edge_external} tell us that in order to obtain a valid reconstruction we have to attach DoFs of vector fields also to geometric elements defined by $\bm s_{n,e}$ and $\bm l_{n,f}$.
In the numerical examples we show how the proposed sparse inverse mass matrices are used in formulations to solve boundary value problems.
In this case, the additional DoFs on the boundary are either known because of boundary conditions or not used in the laws enforced by the formulation.
We emphasize that the additional DoFs on the boundary appears as a canonical and necessary construction to fulfil the geometrical formulas \cref{first_geom} and \cref{second_geom}.

\subsection{Handling material parameter discontinuities inside dual cells}
In this section we are interested in the discretization of problems with discontinuous material parameters, which may appear for example when inhomogeneous Poisson or wave propagation problems are considered.
In order to apply the setting detailed in \cref{dual_mass_matrix}, a key requirement is that the material parameters are uniform inside dual cells. In this case, discontinuities of the material parameters do not create any complication if the grid is chosen in such a way that the interfaces between different materials are aligned with the boundaries of the dual cells. Hence, in this case, the assignment of the material parameters must be based on the dual cells. This is a common practice in the Finite Volume (FE) literature \cite{Eymard2000}.
We observe that such an assignment is correctly defined since dual cells provide a partition of the computational domain.

However, in most cases, it is desired that the material parameters are constant on cells of the primal grid and that they can be arbitrary discontinuous across interfaces between cells.
To each $\C \in C$ is assigned a possibly different material parameter $\wrestr{\material}{\C}$.
As a result, the material parameters of each dual cell are no longer uniform.
In order to handle this kind of situations we can use two different strategies which differ on how the constant vector field is reconstructed inside dual cells.
The first strategy is to use a \emph{weighted average} of the material property in the scheme. The reconstructed vector field is uniform and the construction technique proposed in \cref{dual_mass_matrix} can be applied, by using a suitable choice of the material parameter.
The second strategy consists of using a \emph{piecewise constant vector field} representation inside dual cells, which accounts for discontinuities of the vector field components due to material parameters discontinuities.
\subsubsection{Weighted average}
\label{weighted_material}
Let us consider a dual cell $\DC_n \in \dual C$.
To motivate the definition of the new material parameter, let us consider a uniform vector field $\bm u$ defined in $\DC$. Then, we define the vector field $\bm v$ whose restriction to each cell $\C \in C(n)$ satisfy $\wrestr{\bm v}{\C} \coloneqq \wrestr{\material}{\C} \wrestr{\bm u}{\C}$. Thus, $\bm v$ is a piecewise constant vector field that is spatially constant in each $\C \in C(n)$.
The weighted average of $\bm v$ over $\DC$ is
\begin{align}
\bar{\bm v}
=\frac{1}{\abs{\DC}} \sum_{\C \in C(n)} \abs{\wrestr{\DC}{\C}} \, \wrestr{\material}{\C} \bm u
= \Bigg(\frac{1}{\abs{\DC}} \sum_{\C \in C(n)} \abs{\wrestr{\DC}{\C}} \, \wrestr{\material}{\C} \Bigg) \, \bm u ,
\end{align}
where we have used the fact that $\bm u$ is uniform in $\DC$.
To each dual cell $\DC$ we assign a new material parameter defined by a weighted average as follows
\begin{equation}
\wrestr{\dual {\material}}{\DC} := \Big(\frac{1}{\abs{\DC}} \sum_{\C \in C(n)} \abs{\wrestr{\DC}{\C}} \, \wrestr{\material}{\C} \Big)^{-1}.
\label{material_parameter}
\end{equation}
We use \cref{material_parameter} as the material parameter appearing in the expressions of the local mass matrices in \cref{local_mass_dual}.
As a consequence, if the material parameter is constant over the whole computational domain, the same material parameter is assigned to every dual cell.
\subsubsection{Piecewise constant vector field}
\label{code}
Let us consider a dual cell $\DC_n \in \dual C$.
Since material parameters differ on each $\C \in C$, we use a piecewise constant vector field representation that is spatially constant on each $\C \in C(n)$.
Thus, let us consider two pairs of piecewise constant vector fields defined on $\DC$, namely $(\bm u, \bm v)$ and $(\bm w, \bm z)$.
The restriction of the vector fields $\bm u, \bm v, \bm w, \bm z$ to a cell $\C \in C(n)$ satisfy the following constitutive relations
\begin{align}
\wrestr{\bm v}{\C} = \wrestrs{\material}{1}{\C} \wrestr{\bm u}{\C},\label{dual1cody}
\\
\wrestr{\bm z}{\C} = \wrestrs{\material}{2}{\C} \wrestr{\bm w}{\C},\label{dual2cody}
\end{align}
where $\wrestrs{\material}{1}{\C},\wrestrs{\material}{2}{\C}$ are two symmetric positive definite matrices of order $3$.
In addition, let us suppose that $\bm u, \bm z$ and $\bm v, \bm w$ preserve the \emph{normal} and \emph{tangential} components over the boundaries of $\wrestr{\DC}{\C}$, respectively.
Thus, the restriction of DoFs associated with the vector fields defined in $\DC$ are well defined.
In particular, we attach DoFs to geometric elements of the primal and dual grid as follows $\dofFD{\bm u} = \mFD \bm u$, $\dofDED{\bm v} = \mDED \bm v$, $\dofED{\bm w} = \mED \bm w$, $\dofDFD{\bm z} = \mDFD \bm z$.

The idea underlying the construction is to mimic the reasoning proposed in \cref{local_mass_matrices} and \cref{dual_mass_matrix} to construct global mass matrices starting from local mass matrices defined on cells and dual cells, respectively.
To be more precise, "local" contributions from matrices $\MFR,\MER$ restricted to each region $\wrestr{\DC}{\C}$ are assembled to construct "global" mass matrices $\MFD,\MED$ on the dual cell $\DC$.
The analogy is well defined since the regions $\wrestr{\DC}{\C}$ provide a partition of the cell $\DC$ and the material is uniform on each of them. This latter assumption, allows us to apply the construction of a consistent term using the recipe detailed in \cref{dual_mass_matrix} to each region $\wrestr{\DC}{\C}$.



Thus, by applying \cref{reconsruction_formula_tensor}, a \emph{consistent} mass matrix $\MER$ restricted to $\wrestr{\DC}{\C}$ is given by
\begin{align}
\MER
= \frac{(\mDFDR+\mSR) \wrestr{\material}{\C} (\mDFDR+\mSR)^T}{\abs{\restr{\DC}{\C}}}.
\label{first_new}
\end{align}

A local mass matrix $\MED$ is obtained by assembling the local contributions of each $\MER$.
$\MED$ constructed in this way is symmetric and positive definite since each matrix \cref{first_new} has rank $3$.
Moreover, $\MED$ satisfy the following \emph{piecewise consistency property}, that is
\begin{align}
\dofDFD{\bm z} = \MED \dofED{\bm w}
\end{align}
holds exactly for any pair of piecewise constant vector fields $\bm w,\bm z$ satisfying \cref{dual2cody}.
This follows from the fact that the normal component of $\bm z$ is preserved over the boundaries of $\wrestr{\DC}{\C}$ and hence the geometric boundary terms cancel out, as shown in the proof of \cref{reconstruction_dual_simplicial}.
We observe that the piecewise consistency property implies the consistency property \cref{cons1}.
A local mass matrix $\MDF$ is obtained by computing the algebraic inverse of $\MED$,
\begin{align}
\MDF = {\MED}^{-1}.
\end{align}

In a similar way, by applying \cref{reconsruction_formula_tensor}, a \emph{consistent} mass matrix $\MFR$ restricted to $\wrestr{\DC}{\C}$ is given by
\begin{align}
\MFR
= \frac{(\mDEDR+\mLR) \wrestr{\material}{\C} (\mDEDR+\mLR)^T}{\abs{\restr{\DC}{\C}}}.
\label{second_new}
\end{align}

A local mass matrix $\MFD$ is obtained by assembling the local contributions of each $\MFR$.
$\MFD$ constructed in this way is symmetric and positive definite since each matrix \cref{second_new} has rank $3$.
Moreover, $\MFD$ satisfy the following \emph{piecewise consistency property}, that is
\begin{align}
\dofDED{\bm v} = \MFD \dofFD{\bm u}
\end{align}
holds exactly for any pair of piecewise constant vector fields $\bm v,\bm u$ satisfying \cref{dual1cody}.
This follows from the fact that the tangential component of $\bm v$ is preserved over the boundaries of $\wrestr{\DC}{\C}$ and hence the geometric boundary terms cancel out, as shown in the proof of \cref{reconstruction_dual_simplicial}.
We observe that the piecewise consistency property implies the consistency property \cref{cons2}.
Then, a local mass matrix $\MDF$ is obtained by computing the algebraic inverse of $\MED$,
\begin{align}
\MDE = {\MFD}^{-1}.
\end{align}

We observe that even if the material parameter is uniform inside $\DC$, using this approach it is in any case necessary to compute the inverses of matrices $\MED,\MFD$.
Instead, in \cref{dual_mass_matrix}, explicit expressions of $\MDF,\MDE$ are derived.
\subsubsection{Hybrid approach to handle material discontinuities inside dual cells}
\label{hybrid}
The proposed approach to handle discontinuities of the material parameter inside dual cells is the following. First, if material parameters are uniform inside a dual cell we can apply the construction detailed in \cref{dual_mass_matrix}.
Otherwise, we apply one among the approaches described in \cref{weighted_material}, \cref{code}.
We point out that these constructions are applied only to cells which lie at the intersection between regions enclosing different materials.
Thus, the required computational effort to compute local inverses in the approach described in \cref{code} is negligible in practice.

An important remark is that if we want to satisfy exactly a \emph{multi-material patch test} (see \cref{numerical_results}), it is necessary in our hybrid approach to resort to the construction proposed in \cref{code}.
This is because, being the reconstructed field piecewise constant, it can represent exactly the discontinuities of the vector field components that appear at interfaces between different materials.
Instead, in \cref{weighted_material}, the reconstructed field is constant.

\section{Numerical results}
\label{numerical_results}
Beside all the other applications of sparse inverse mass matrices, to validate the construction proposed in this paper we concentrate on the solution of Poisson boundary value problems formulated with one unknown per element.

We first verified that the proposed technique is able to pass a patch test on grids made by general tetrahedra.
Then, a stationary current conduction problem representing a square resistor is solved. The stationary current conduction is a Poisson problem in region $\Omega$ of the 3-D Euclidean space
\begin{subequations}
\begin{align}
&\curl \bm E = \bm 0,\label{faraday} \\
&\div \bm J = 0,\label{solenoidal} \\
&\bm J = \sigma \bm E, \label{constitutive}
\end{align}
\end{subequations}
where $\sigma$ is the material parameter electric conductivity, $\bm E$ and $\bm J$ are the \emph{conservative} electrostatic field and the current density vectors, respectively.
The electric conductivity $\sigma$ is assumed to be a positive scalar value which is piecewise uniform in each material region.
The region boundary $\partial \Omega$ is partitioned into a set of $N^{\mathrm{i}}$ surfaces of perfect insulators $\partial \Omega^{\mathrm{i}}_k$, and  a set of $N^{\mathrm{c}}+1$ disjoint equipotential surfaces of perfect conductors $\partial \Omega^{\mathrm{c}}_k$ (usually called electrodes):
\begin{align}
\label{m4}
 \partial \Omega=
\sum_{k=1}^{N^{\mathrm{i}}} \partial \Omega^{\mathrm{i}}_k +
\sum_{k=0}^{N^{\mathrm{c}}} \partial \Omega^{\mathrm{c}}_k.
\end{align}
Electrode $\partial \Omega^{\mathrm{c}}_0$ is considered as reference for all the voltages of the remaining electrodes, that are supposed to be assigned. $\bm{J} \cdot {\bm n}=0$ is set as boundary conditions (b.c.) on each $\partial \Omega^{\mathrm{i}}_k$, where ${\bm n}$ is the outwards oriented normal unit vector of $\partial \Omega$.

Since $\curl \bm E = \bm 0$ and the electrostatic field is a conservative field, we introduce the scalar potential $U$ such that $\bm E = - \grad U$. Similarly, since $\div \bm J = \bm 0$, we introduce a vector potential $\bm T$  such that $\bm J = \curl \bm T$.

\subsection{Survey of standard formulations}
There are several formulations to numerically solve Poisson problems like \eqref{faraday}-\eqref{solenoidal}-\eqref{constitutive}.
The most common one is the FEM formulation based on the scalar potential $SP$ expanded with the classical nodal basis functions. With this formulation the unknowns $\GdofN{U}$, are the DoFs of the scalar potentials $U$ sampled on the grid nodes. The voltages $\GdofE{\bm E} = -\mat G \GdofN{U}$, are associated with the grid edges, where $\mat G$ is the edge-node incidence matrix. Finally, the DoFs of $\bm J$ are attached to dual faces as \cite{howweak} shows. Faraday's law \eqref{faraday} is enforced implicitly by the scalar potential, whereas the solenoidality of the current \eqref{solenoidal} is enforced on the boundary of each dual cell by the linear system, see for example \cite{ijnme}.

Other possibilities emerge as we exchange the association of physical variables between the primal and dual grids. In \textit{complementary} formulations the DoFs of $\bm J$ are attached to faces of the primal grid, whereas the scalar potential DoFs are attached to dual nodes. Consequently, the voltages $\GdofDE{\bm E}=-\dualG \GdofDN{U}$, are associated with dual edges.
There are two methods to obtain a complementary formulation.
The first \textit{complementary} formulation $VP$  fulfill \eqref{solenoidal} by using the electric vector potential $\mathbf{T}$ \cite{ijnme} such that $\GdofF{\bm J}= \mat C \GdofE{\bm T}$ (We remark that in general the so-called \emph{relative cohomology theory} is required to present this formulation, see \cite{ijnme}).
We note that in the $VP$ formulation the role of physical laws are exchanged w.r.t. the $SP$ formulation since Faraday's law \cref{faraday} is enforced with a linear system in the VP formulation.

There is a second method to produce complementary formulations that we call \emph{complementary-dual}: they are complementary formulation but still use the scalar potential $\GdofDN{U}$, which is sampled on grid dual nodes, one-to-one with grid cells. An effective complementary-dual formulation, which is algebraically equivalent to the $VP$ formulation, is the mixed-hybrid $MH$ formulation \cite{brezzicmame2007}, \cite{mh}.

\subsection{Formulations with one unknown per element}
There are other complementary-dual formulations, much less explored in the literature, that feature one unknown per element \cite{Eymard2000}, \cite{vohralik}, \cite{bellina}. As pointed out in \cite{vohralik}, there has been a long-standing interest to reduce the system to one potential value per element, to reduce unknowns and obtain a positive-definite system.

In \cite{vohralik} two approaches are proposed. The first one requires the solution of local problems on patches of elements to produce local flux expressions. The second one is based on the use of the algebraic inverse of the mass matrix and has been deemed as impossible in \cite{vohralik} because the ``the inverse of a mass matrix is not sparse''.
We follow the latter approach enabled by the efficient construction of sparse inverse of mass matrices proposed for the first time in this paper. To show the details, let us consider the dual scalar potential formulations $DSP$ obtained by writing the problem in the geometric framework \cite{howweak} as
%
\begin{align}
  &\mat C^T \GdofDE{\bm E}=\bm 0\label{m1d}\\
  &\mat D \GdofF{\bm J} =\bm 0 \label{m2d}\\
  &\GdofF{\bm J} = \GMDE \GdofDE{\bm E},\label{m3d}
\end{align}
%
where $\mat D$ is the cell-face incidence matrix, $\mat C$ is the face-edge incidence matrix and $\GMDE$ is the inverse mass matrix of $\GME$, which maps DoFs of $\bm E$ attached to dual edges to DoFs of $\bm J$ attached to primal faces. Thus, we can interpret $\GMDE$ as a \emph{dual conductance matrix}.
To implicitly satisfy \eqref{m1d}, the scalar potential $\GdofDN{U}$ in the dual nodes is introduced through
\begin{equation}\label{pot22}
\GdofDE{\bm E}
=-\dualG \GdofDN{U} + \GdofDE{\bm E}_s,
\end{equation}
%
where $\GdofDE{\bm E}_s$ is introduced to take into account Dirichlet b.c., so that
$\dualC \GdofDE{\bm E}_s=\bm 0$ \cite{bossavitmh}, \cite{mh}. Its construction is straightforward and detailed in \cite{mh}. By substituting \eqref{pot22} and \eqref{m3d} into \eqref{m2d}, one gets
\begin{equation}
\label{C2019_eq1}
(\mat D \GMDE \mat D^T)\GdofDN{U}
=\mat D \GMDE \GdofDE{\bm E}_s,
\end{equation}
having the scalar potential on dual nodes as unknowns.

A symmetric and positive-definite $\GMDE$ may be computed as $\GMDE= {\GME}^{-1}$,
where the Raviart--Thomas mass matrix $\GME$ can be interpreted as a \emph{resistance mass matrix}.
Yet, the obtained matrix $\GMDE$ using this recipe is fully populated, hence not usable in practice \cite{vohralik}.
Another solution is enabled by the geometric construction of $\GMDE$ by using \cref{facedual}, which produces a symmetric and positive definite matrix for any tetrahedral grid given as input.
We point out that the additional voltages DoFs on the boundary are zero on the electrodes, because electrodes are equipotential by hypothesis. The other additional voltages DoFs on the boundary are not needed because the boundary conditions are applied on the currents of the corresponding boundary faces.

We remark that another formulation arises that we may call dual vector potential $DVP$, which is the dual of the $VP$ formulation. I.e. one obtains a system in which the unknowns are the DoFs of the vector potential attached on dual edges. This formulation is not presented in detail since it has been verified that it is far to be competitive with the others in terms of computational efficiency.

\subsection{Classical patch tests}
\begin{figure}[!h]
	\centering
	\includegraphics[scale=0.5]{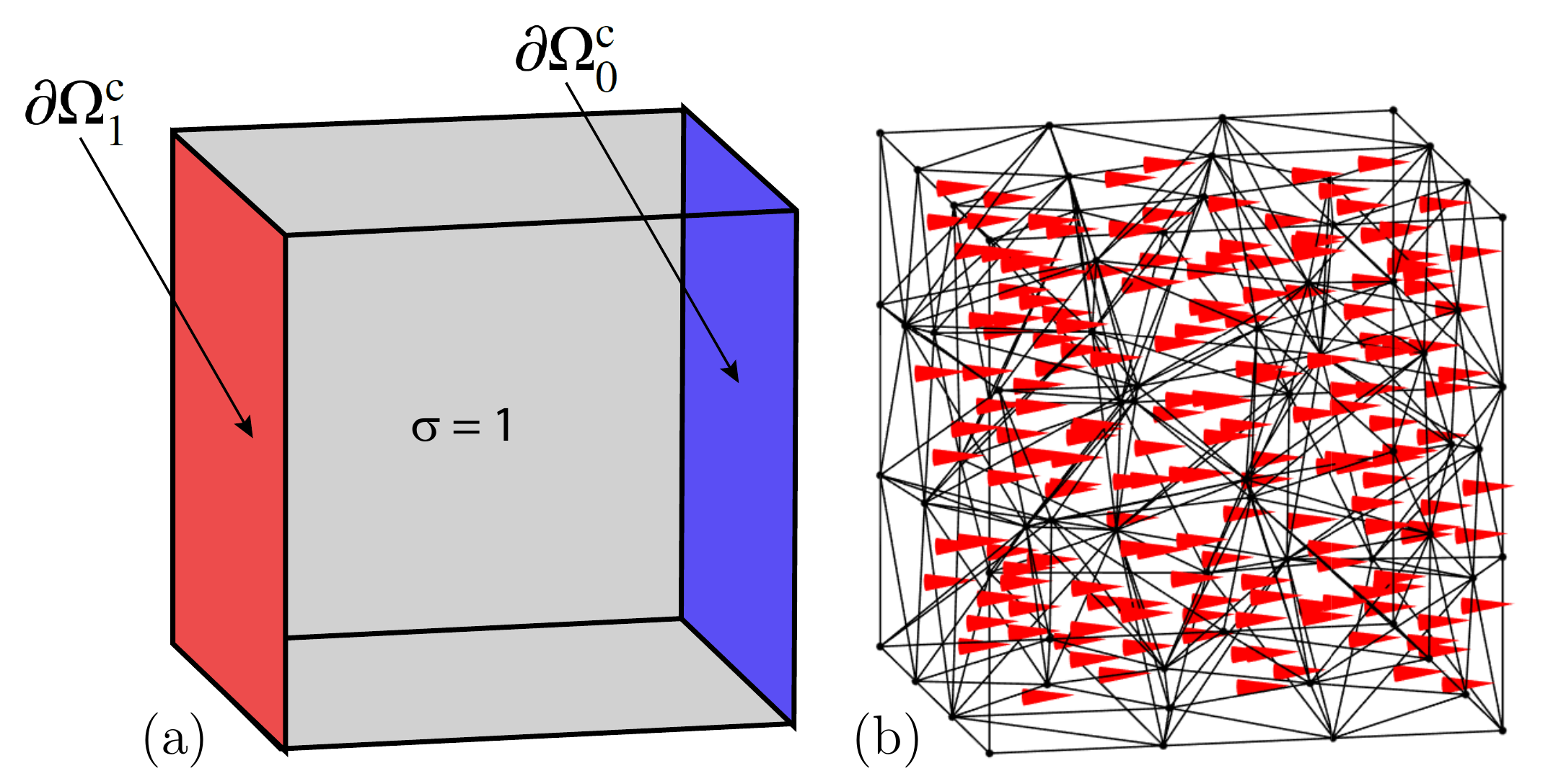}
	\caption{(a) The definition of a patch test as the solution of an electric conduction problem inside a planar resistor. (b) All mentioned formulations are able to retrieve the analytical solution up to machine precision or iterative solver tolerance.}
	\label{patch}
\end{figure}

The patch tests are Poisson problems designed in such a way that their analytical solution is uniform. By interpreting the Poisson problems as direct current conduction problems, a simple way to produce a patch test is to consider a planar resistor, as described in Fig. \ref{patch}(a). It has been verified that the $SP$, $VP=MH$, and $DSP$ formulations produce the analytical solution as represented in \cref{patch}(b).

\subsection{Multi-material patch tests}
\begin{figure}[!h]
	\centering
	\includegraphics[scale=0.35]{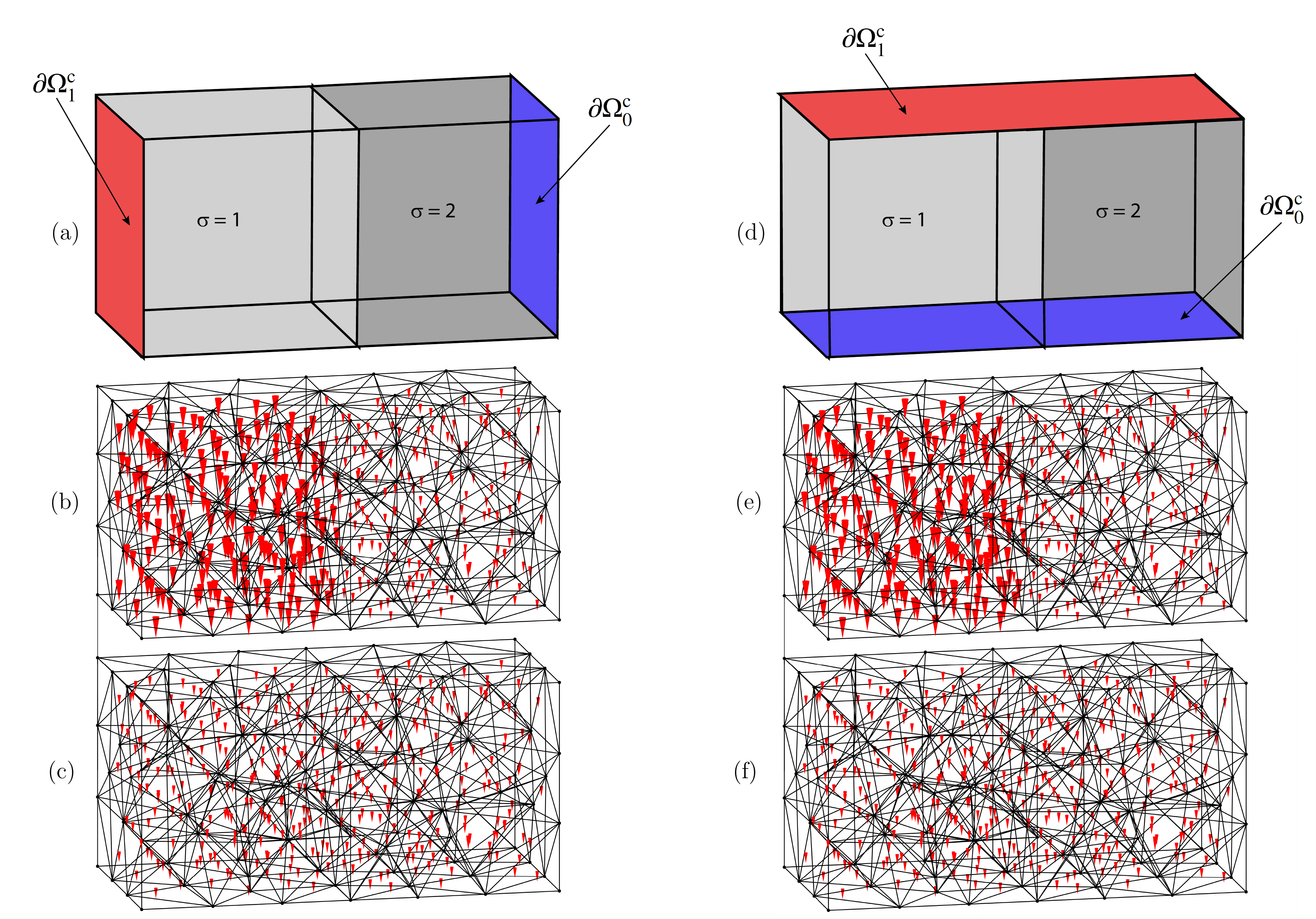}
	\caption{(a) The \emph{series} multi-material patch test. (b) Electrostatic field $\bm E$ obtained with the series. (c) Current density field $\bm J$ obtained with the series. (d) The \emph{parallel} multi-material patch test. (e) Electric field $\bm{E}$ obtained with the parallel. (f) Current density field $\bm{J}$ obtained with the parallel.}
	\label{mpatch}
\end{figure}

The multi-material patch tests are direct current conduction problems designed in such a way that their solution is \emph{piecewise uniform}. A simple way to produce multi-material patch test is to consider a resistor with two conductors with different material properties placed in \emph{series} or in \emph{parallel} as described in \cref{mpatch}(a) and \cref{mpatch}(d). A voltage of $u=1\,$V is enforced between the reference electrode $\partial \Omega^{\mathrm{c}}_0$ and the electrode $\partial \Omega^{\mathrm{c}}_1$.

In the first multi-material patch test, two different materials with different material properties are placed in \emph{series}. From the result represented in \cref{mpatch}(b-c), we conclude that the tangential component of the electric field $\bm{E}$ is conserved across the material interface, whereas the tangential component of the current density field $\bm{J}$ jumps. This result holds for all formulations.

In the second multi-material patch test, two different materials with different material properties are placed in \emph{parallel}. From the result represented in \cref{mpatch}(e-f), we conclude that the normal component of the current density $\bm{J}$ is conserved across the material interface, whereas the normal component of the electric field $\bm{E}$ jumps. This result holds for all formulations.

\subsection{Square resistor benchmark}
\begin{figure}[!t]
	\centering
	\includegraphics[scale=0.5]{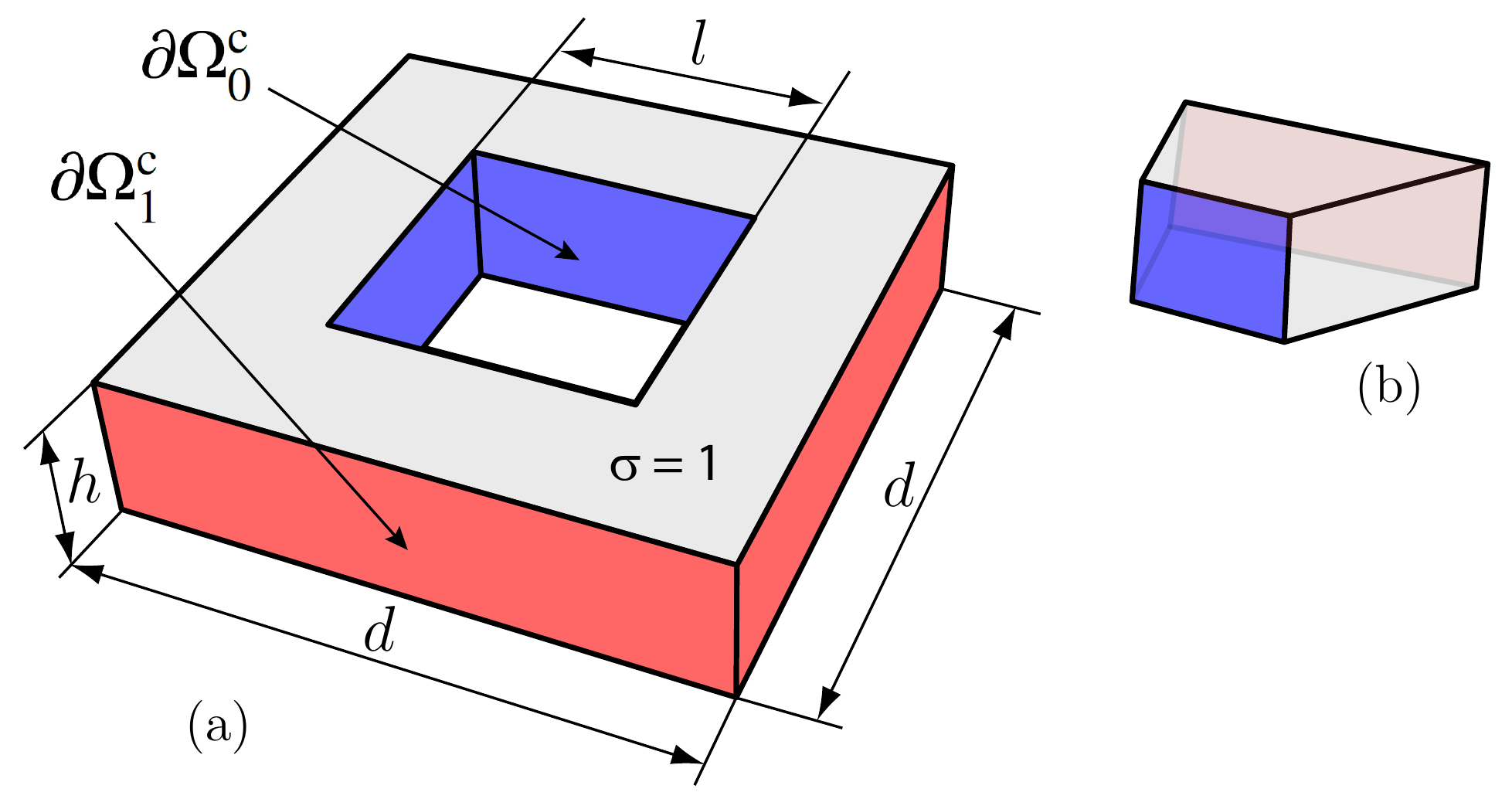}
	\caption{(a) The geometry of the square resistor benchmark ($h = 1\,$m, $d = 4\,$m, and
$l = 2\,$m). The two electrodes are depicted in red and blue. (b) Thanks to the symmetry, the computational domain has been reduced to one in eight of the resistor.}
	\label{sres}
\end{figure}

As a more complicated example, we compute the conductance $G$ of a square resistor \cite{mh}, see \cref{sres}(a). A voltage of $u=1\,$V is enforced between the reference electrode $\partial \Omega^{\mathrm{c}}_0$ and the electrode $\partial \Omega^{\mathrm{c}}_1$, which are placed in the two lateral surfaces of the solid square torus. The conductor placed inside the solid torus has an electrical conductivity $\sigma=1\,$S/m. This problem, like most industrial problems, exhibits a singular analytical solution. Yet, the analytical solution is available ($G=10.23409256\,$S).
The conductance $G$ is extracted by computing the total dissipated power $P=G\,u^2=\int_{\Omega} \frac{|\bm{J}|^2}{\sigma}$ or alternatively by using the Ohm's law $G=i/u$, where $i$ is the current that flows between the two electrodes.

\begin{figure}[!t]
	\centering
	\includegraphics[scale=0.7]{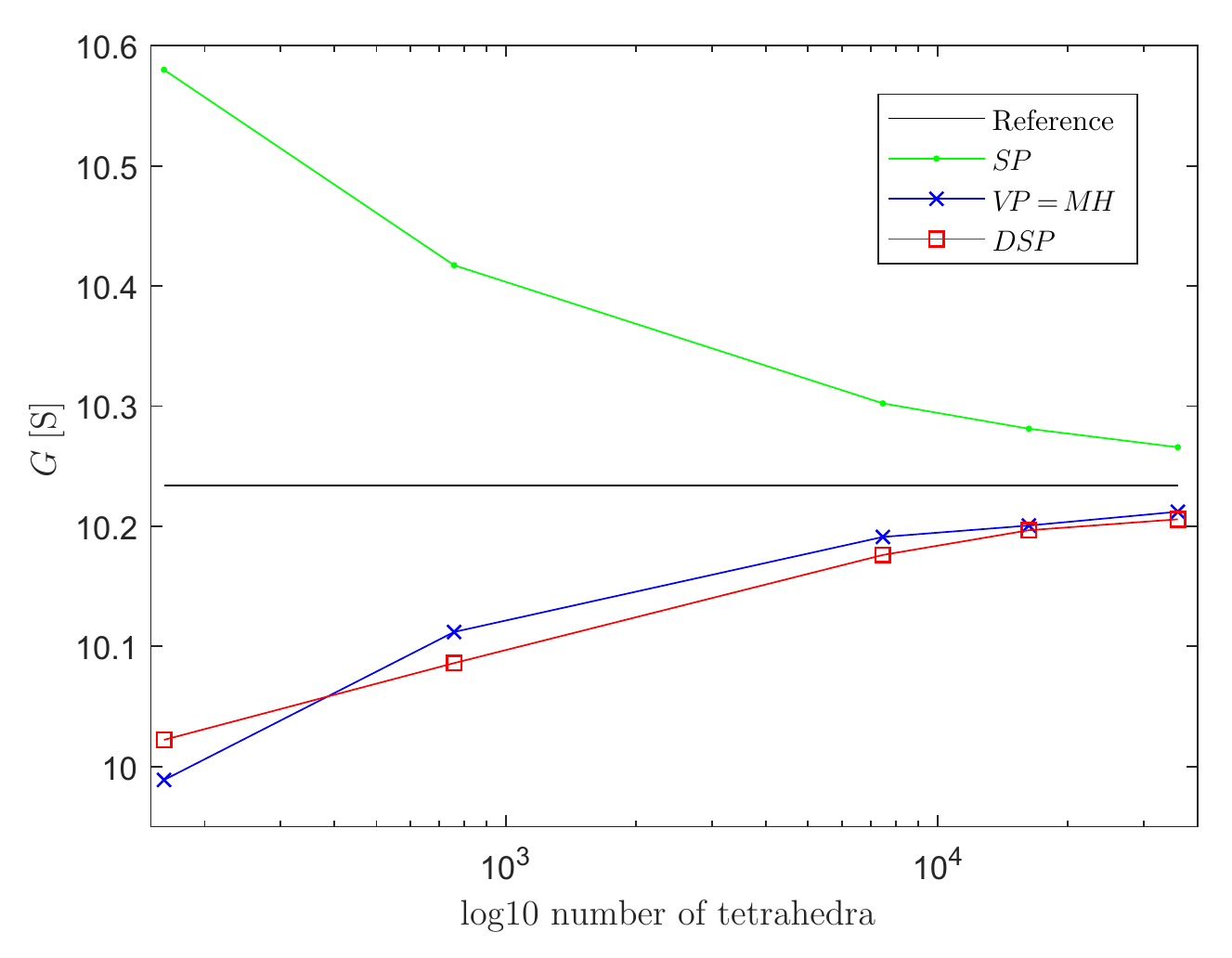}
	\caption{Results for the square resistor benchmark.}
	\label{combres}
\end{figure}

We recall that the $VP$ and $MH$ formulations produce the same results given that they are algebraically equivalent \cite{mh}.
The conductance of the square resistor has been evaluated by the $SP$, $VP=MH$, and $DSP$ formulations on refined grids. The results are collected in \cref{combres}.
First, it is interesting to note that the results relative to the scalar potential $SP$ and to the vector potential $VP$ or mixed-hybrid $MH$ formulations provide, respectively, the upper and lower bounds for the conductance \cite{synge}, \cite{ijnme}.
This property is called \emph{complementarity} in the computational electromagnetics community \cite{bossavitbook}, \cite{mh}.


\section{Conclusions}
\label{conclusions}
In this paper, we have proposed a new construction of sparse inverse edge and face mass matrices for arbitrary tetrahedral grids.
The proposed framework unifies the construction of both mass matrices and their inverses and is based on novel geometric reconstruction formulas, from which, local mass matrices are defined as the sum of a consistent and a stabilization term.
The consistent term is defined geometrically and explicitly, thus providing a sensible speedup and an easier implementation, whereas the stabilization term is constructed according to a well-established design strategy.
The resulting expressions for both types of local mass matrices are highly symmetric, thus highlighting the duality relationship between the pair of grids.
Global mass matrices are then constructed using a standard assembling process.
In this way we obtain a sparse matrix representation also for the inverse mass matrix.
This avoids explicitly computing the algebraic inverse of the mass matrix which is found to be typically dense.
A key aspect of the proposed method is that being based on geometric elements of the barycentric dual grid, it can be applied to arbitrary tetrahedral grids.
Two different strategies are proposed to deal with the case of problems involving discontinuous material parameters for the case where the interfaces are aligned with the tetrahedral grid, and an hybrid approach between the two is employed in practical applications.
We tested the newly derived inverse mass matrices by solving a Poisson problem and we verified that patch test is passed, even when material discontinuities are present.
A real case problem involving a square resistor benchmark is analyzed.
Given that reconstruction formulas in each dual cell provide a first-order approximation of a vector field, we expect that first-order error estimates can be derived for the numerical scheme.
%
%
%
%

\bibliographystyle{elsarticle-num}
\bibliography{ref}

\begin{thebibliography}{10}
\expandafter\ifx\csname url\endcsname\relax
  \def\url#1{\texttt{#1}}\fi
\expandafter\ifx\csname urlprefix\endcsname\relax\def\urlprefix{URL }\fi
\expandafter\ifx\csname href\endcsname\relax
  \def\href#1#2{#2} \def\path#1{#1}\fi

\bibitem{tontiold}
E.~Tonti, On the formal structure of physical theories, Monograph of the
  Italian National Research Council, 1975.

\bibitem{Yee1966}
K.~S. Yee, {Numerical Solution of Initial Boundary Value Problems Involving
  Maxwell's Equations in Isotropic Media}, IEEE Transactions on Antennas and
  Propagation 14~(3) (1966) 302--307.

\bibitem{Tonti2001}
E.~Tonti, {A direct discrete formulation of field laws: The cell method}, CMES
  - Computer Modeling in Engineering and Sciences 2~(2) (2001) 237--258.
\newblock \href {https://doi.org/10.3970/cmes.2001.002.237}
  {\path{doi:10.3970/cmes.2001.002.237}}.

\bibitem{Clemens2001}
M.~Clemens, T.~Weiland, {Discrete electromagnetism with the finite integration
  technique}, Progress in Electromagnetics Research 32 (2001) 65--87.
\newblock \href {https://doi.org/10.2528/PIER00080103}
  {\path{doi:10.2528/PIER00080103}}.

\bibitem{cmame}
L.~Codecasa, R.~Specogna, F.~Trevisan,
  \href{http://www.sciencedirect.com/science/article/pii/S0045782508004234}{Base
  functions and discrete constitutive relations for staggered polyhedral
  grids}, Comput. Meth. Appl. Mech. Eng. 198~(9) (2009) 1117 -- 1123.
\newblock \href {https://doi.org/https://doi.org/10.1016/j.cma.2008.11.021}
  {\path{doi:https://doi.org/10.1016/j.cma.2008.11.021}}.
\newline\urlprefix\url{http://www.sciencedirect.com/science/article/pii/S0045782508004234}

\bibitem{jcp}
L.~Codecasa, R.~Specogna, F.~Trevisan,
  \href{http://www.sciencedirect.com/science/article/pii/S0021999110003384}{A
  new set of basis functions for the discrete geometric approach}, J. Comput.
  Phys. 229~(19) (2010) 7401 -- 7410.
\newblock \href {https://doi.org/https://doi.org/10.1016/j.jcp.2010.06.023}
  {\path{doi:https://doi.org/10.1016/j.jcp.2010.06.023}}.
\newline\urlprefix\url{http://www.sciencedirect.com/science/article/pii/S0021999110003384}

\bibitem{bonellem2an}
J.~Bonelle, A.~Ern,
  \href{http://www.numdam.org/item/M2AN_2014__48_2_553_0}{Analysis of
  compatible discrete operator schemes for elliptic problems on polyhedral
  meshes}, ESAIM: Mathematical Modelling and Numerical Analysis -
  Mod\'elisation Math\'ematique et Analyse Num\'erique 48~(2) (2014) 553--581.
\newblock \href {https://doi.org/10.1051/m2an/2013104}
  {\path{doi:10.1051/m2an/2013104}}.
\newline\urlprefix\url{http://www.numdam.org/item/M2AN_2014__48_2_553_0}

\bibitem{Lipnikov2014}
K.~Lipnikov, G.~Manzini, M.~Shashkov, {Mimetic finite difference method}, J.
  Comput. Phys. 257~(PB) (2014) 1163--1227.
\newblock \href {https://doi.org/10.1016/j.jcp.2013.07.031}
  {\path{doi:10.1016/j.jcp.2013.07.031}}.

\bibitem{Bossavit1988}
A.~Bossavit, {Whitney forms: A class of finite elements for three-dimensional
  computations in electromagnetism}, IEE Proceedings A: Physical Science.
  Measurement and Instrumentation. Management and Education. Reviews 135 pt
  A~(8) (1988) 493--500.
\newblock \href {https://doi.org/10.1049/ip-a-1.1988.0077}
  {\path{doi:10.1049/ip-a-1.1988.0077}}.

\bibitem{bossavitbook}
A.~Bossavit, Computational Electromagnetism, Academic Press, 1998.

\bibitem{howweak}
A.~{Bossavit}, How weak is the "weak solution" in finite element methods?, IEEE
  Transactions on Magnetics 34~(5) (1998) 2429--2432.

\bibitem{tontinew}
E.~Tonti, The mathematical structure of classical and relativistic physics. A
  general classification diagram, Birkh\"{a}user, Basel, 2013.

\bibitem{kettunen}
T.~{Tarhasaari}, L.~{Kettunen}, A.~{Bossavit}, Some realizations of a discrete
  hodge operator: a reinterpretation of finite element techniques [for em field
  analysis], IEEE Trans. Magn. 35~(3) (1999) 1494--1497.

\bibitem{Codecasa2008}
L.~Codecasa, M.~Politi, {Explicit, consistent, and conditionally stable
  extension of FD-TD to tetrahedral grids by FIT}, IEEE Transactions on
  Magnetics 44~(6) (2008) 1258--1261.
\newblock \href {https://doi.org/10.1109/TMAG.2007.916310}
  {\path{doi:10.1109/TMAG.2007.916310}}.

\bibitem{teixeira}
B.~{He}, F.~L. {Teixeira}, Differential forms, galerkin duality, and sparse
  inverse approximations in finite element solutions of maxwell equations, IEEE
  Transactions on Antennas and Propagation 55~(5) (2007) 1359--1368.

\bibitem{Codecasa2018}
L.~Codecasa, B.~Kapidani, R.~Specogna, F.~Trevisan, {Novel FDTD technique over
  tetrahedral grids for conductive media}, IEEE Transactions on Antennas and
  Propagation 66~(10) (2018) 5387--5396.
\newblock \href {https://doi.org/10.1109/TAP.2018.2862244}
  {\path{doi:10.1109/TAP.2018.2862244}}.

\bibitem{bell}
N.~Bell, L.~N. Olson,
  \href{https://onlinelibrary.wiley.com/doi/abs/10.1002/nla.577}{Algebraic
  multigrid for k-form laplacians}, Numerical Linear Algebra with Applications
  15~(2‐3) (2008) 165--185.
\newblock \href
  {http://arxiv.org/abs/https://onlinelibrary.wiley.com/doi/pdf/10.1002/nla.577}
  {\path{arXiv:https://onlinelibrary.wiley.com/doi/pdf/10.1002/nla.577}}, \href
  {https://doi.org/10.1002/nla.577} {\path{doi:10.1002/nla.577}}.
\newline\urlprefix\url{https://onlinelibrary.wiley.com/doi/abs/10.1002/nla.577}

\bibitem{caltagirone}
A.~Lemoine, J.~Caltagirone, M.~Aza{\"{\i}}ez, S.~Vincent,
  \href{https://doi.org/10.1007/s10915-014-9952-8}{Discrete helmholtz-hodge
  decomposition on polyhedral meshes using compatible discrete operators}, J.
  Sci. Comput. 65~(1) (2015) 34--53.
\newblock \href {https://doi.org/10.1007/s10915-014-9952-8}
  {\path{doi:10.1007/s10915-014-9952-8}}.
\newline\urlprefix\url{https://doi.org/10.1007/s10915-014-9952-8}

\bibitem{Hirani2003}
A.~N. Hirani, {Discrete Exterior Calculus Thesis by}, Tech. rep. (2003).

\bibitem{Eymard2000}
R.~Eymard, T.~Gallou{\"{e}}t, R.~Herbin, andRaph{\`{a}}~ele Herbin,
  \href{https://hal.archives-ouvertes.fr/hal-02100732v2}{{Handbook of Numerical
  Analysis, 9780444503503. which appeared in Handbook of Numerical Analysis}},
  P.G. Ciarlet, J.L. Lions eds 7 (2000) 713--1020.
\newblock \href {https://doi.org/10.1016/S1570-8659(00)07005}
  {\path{doi:10.1016/S1570-8659(00)07005}}.
\newline\urlprefix\url{https://hal.archives-ouvertes.fr/hal-02100732v2}

\bibitem{bajaj}
A.~Gillette, C.~Bajaj,
  \href{http://www.sciencedirect.com/science/article/pii/S001044851100159X}{Dual
  formulations of mixed finite element methods with applications},
  Computer-Aided Design 43~(10) (2011) 1213 -- 1221, solid and Physical
  Modeling 2010.
\newblock \href {https://doi.org/https://doi.org/10.1016/j.cad.2011.06.017}
  {\path{doi:https://doi.org/10.1016/j.cad.2011.06.017}}.
\newline\urlprefix\url{http://www.sciencedirect.com/science/article/pii/S001044851100159X}

\bibitem{bellina}
F.~{Bellina}, R.~{Specogna}, Diagonal material matrices for arbitrary
  simplicial meshes for solving poisson problems with one unknown per element,
  IEEE Transactions on Magnetics 56~(2) (2020) 1--4.

\bibitem{bossavitmh}
A.~{Bossavit}, Mixed-hybrid methods in magnetostatics: complementarity in one
  stroke, IEEE Transactions on Magnetics 39~(3) (2003) 1099--1102.

\bibitem{hirani2}
M.~S. Mohamed, A.~N. Hirani, R.~Samtaney, Numerical convergence of discrete
  exterior calculus on arbitrary surface meshes, International Journal for
  Computational Methods in Engineering Science and Mechanics 19~(3) (2018)
  194--206.

\bibitem{iceaa}
L.~{Codecasa}, R.~{Specogna}, F.~{Trevisani}, The discrete geometric approach
  for wave propagation problems, in: 2009 International Conference on
  Electromagnetics in Advanced Applications, 2009, pp. 59--62.
\newblock \href {https://doi.org/10.1109/ICEAA.2009.5297610}
  {\path{doi:10.1109/ICEAA.2009.5297610}}.

\bibitem{beltman}
R.~Beltman, Mimetic discretizations of the incompressible navier--stokes
  equations for polyhedral meshes, Ph.D. thesis, Department of Mathematics and
  Computer Science, proefschrift. (Nov. 2020).

\bibitem{bonellecad}
J.~Bonelle, D.~A. {Di Pietro}, A.~Ern, {Low-order reconstruction operators on
  polyhedral meshes: Application to compatible discrete operator schemes},
  Computer Aided Geometric Design 35-36 (2015) 27--41.
\newblock \href {https://doi.org/10.1016/j.cagd.2015.03.015}
  {\path{doi:10.1016/j.cagd.2015.03.015}}.

\bibitem{Bochev2007}
P.~B. Bochev, J.~M. Hyman, {Principles of Mimetic Discretizations of
  Differential Operators}, in: Compatible Spatial Discretizations, Springer New
  York, 2007, pp. 89--119.
\newblock \href {https://doi.org/10.1007/0-387-38034-5_5}
  {\path{doi:10.1007/0-387-38034-5_5}}.

\bibitem{munkres}
J.~Munkres, Elements of algebraic topology, Perseus Books, Cambridge, MA, 1984.

\bibitem{Berger}
M.~Berger, M.~Cole, S.~Levy, Geometry I, Universitext, Springer Berlin
  Heidelberg, 2009.

\bibitem{ijnme}
R.~Specogna,
  \href{https://onlinelibrary.wiley.com/doi/abs/10.1002/nme.3089}{Complementary
  geometric formulations for electrostatics}, Int. J. Numer. Meth. Eng. 86~(8)
  (2011) 1041--1068.
\newblock \href
  {http://arxiv.org/abs/https://onlinelibrary.wiley.com/doi/pdf/10.1002/nme.3089}
  {\path{arXiv:https://onlinelibrary.wiley.com/doi/pdf/10.1002/nme.3089}},
  \href {https://doi.org/10.1002/nme.3089} {\path{doi:10.1002/nme.3089}}.
\newline\urlprefix\url{https://onlinelibrary.wiley.com/doi/abs/10.1002/nme.3089}

\bibitem{brezzicmame2007}
F.~Brezzi, K.~Lipnikov, M.~Shashkov, V.~Simoncini,
  \href{http://www.sciencedirect.com/science/article/pii/S0045782507000965}{A
  new discretization methodology for diffusion problems on generalized
  polyhedral meshes}, Comput. Meth. Appl. Mech. Eng. 196~(37) (2007) 3682 --
  3692.
\newblock \href {https://doi.org/https://doi.org/10.1016/j.cma.2006.10.028}
  {\path{doi:https://doi.org/10.1016/j.cma.2006.10.028}}.
\newline\urlprefix\url{http://www.sciencedirect.com/science/article/pii/S0045782507000965}

\bibitem{mh}
R.~{Specogna}, One stroke complementarity for poisson-like problems, IEEE
  Trans. Magn. 51~(3) (2015) 1--4.

\bibitem{vohralik}
M.~VohralÍ\'{\i}k, B.~I. Wohlmuth, Mixed finite element methods:
  implementation with one unknown per element, local flux expressions,
  positivity, polygonal meshes, and relations to other methods, Math. Mod.
  Meth. Appl. S. 23~(05) (2013) 803--838.

\bibitem{synge}
J.~Synge, The hypercircle in mathematical physics, Cambridge University Press,
  Cambridge, 1957.

\end{thebibliography}



\end{document}